\documentclass[reqno,11pt]{amsart}
\usepackage{amsmath}
\usepackage{amsxtra}
\usepackage{amscd}
\usepackage{amsthm}
\usepackage{amsfonts}
\usepackage{amssymb}
\usepackage{eucal}
\usepackage[matrix,arrow,curve]{xy}
\usepackage[dvipsnames]{xcolor}

\theoremstyle{plain}
\newtheorem{Thm}[subsection]{Theorem}
\newtheorem{Cor}[subsection]{Corollary}
\newtheorem{Lem}[subsection]{Lemma}
\newtheorem{Prop}[subsection]{Proposition}
\newtheorem{Conj}[subsection]{Conjecture}

\theoremstyle{definition}
\newtheorem{Def}[subsection]{Definition}

\theoremstyle{remark}

\newtheorem{Rem}[subsection]{Remark}

\newtheorem{conj}{Conjecture}

\numberwithin{equation}{section}
\renewcommand{\rm}{\normalshape}

\newif\ifShowLabels
\ShowLabelstrue
\newdimen\theight
\def\TeXref#1{%
    \leavevmode\vadjust{\setbox0=\hbox{{\tt
        \quad\quad  {\small \rm #1}}}%
    \theight=\ht0
    \advance\theight by \lineskip
    \kern -\theight \vbox to
    \theight{\rightline{\rlap{\box0}}%
    \vss}%
    }}%

\ShowLabelsfalse% comment this out if labels should be printed

\renewcommand{\sec}[2]{\section{#2}\label{S:#1}%
    \ifShowLabels \TeXref{{S:#1}} \fi}
\newcommand{\ssec}[2]{\subsection{#2}\label{SS:#1}%
    \ifShowLabels \TeXref{{SS:#1}} \fi}

\newcommand{\refs}[1]{Section ~\ref{S:#1}}
\newcommand{\refss}[1]{Section ~\ref{SS:#1}}

\newcommand{\reft}[1]{Theorem ~\ref{T:#1}}
\newcommand{\refl}[1]{Lemma ~\ref{L:#1}}
\newcommand{\refp}[1]{Proposition ~\ref{P:#1}}
\newcommand{\refc}[1]{Corollary ~\ref{C:#1}}
\newcommand{\refd}[1]{Definition ~\ref{D:#1}}
\newcommand{\refr}[1]{Remark ~\ref{R:#1}}
\newcommand{\refe}[1]{\eqref{E:#1}}

\newenvironment{thm}[1]%
    { \begin{Thm} \label{T:#1}  \ifShowLabels \TeXref{T:#1} \fi }%
    { \end{Thm} }

\renewcommand{\th}[1]{\begin{thm}{#1} \sl }
\renewcommand{\eth}{\end{thm} }

\newenvironment{lemma}[1]%
    { \begin{Lem} \label{L:#1}  \ifShowLabels \TeXref{L:#1} \fi }%
    { \end{Lem} }
\newcommand{\lem}[1]{\begin{lemma}{#1} \sl}
\newcommand{\elem}{\end{lemma}}

\newenvironment{propos}[1]%
    { \begin{Prop} \label{P:#1}  \ifShowLabels \TeXref{P:#1} \fi }%
    { \end{Prop} }
\newcommand{\prop}[1]{\begin{propos}{#1}\sl }
\newcommand{\eprop}{\end{propos}}

\newenvironment{corol}[1]%
    { \begin{Cor} \label{C:#1}  \ifShowLabels \TeXref{C:#1} \fi }%
    { \end{Cor} }
\newcommand{\cor}[1]{\begin{corol}{#1} \sl }
\newcommand{\ecor}{\end{corol}}

\newenvironment{defeni}[1]%
    { \begin{Def} \label{D:#1}  \ifShowLabels \TeXref{D:#1} \fi }%
    { \end{Def} }
\newcommand{\defe}[1]{\begin{defeni}{#1} \sl }
\newcommand{\edefe}{\end{defeni}}

\newenvironment{remark}[1]%
    { \begin{Rem} \label{R:#1}  \ifShowLabels \TeXref{R:#1} \fi }%
    { \end{Rem} }
\newcommand{\rem}[1]{\begin{remark}{#1}}
\newcommand{\erem}{\end{remark}}

\newenvironment{conjec}[1]%
    { \begin{Conj} \label{Co:#1}  \ifShowLabels \TeXref{Co:#1} \fi }%
    { \end{Conj} }
\renewcommand{\conj}[1]{\begin{conjec}{#1} \sl }
\newcommand{\econj}{\end{conjec}}

\newcommand{\eq}[1]%
    { \ifShowLabels \TeXref{E:#1} \fi
       \begin{equation} \label{E:#1} }
\newcommand{\eeq}{ \end{equation} }

\newcommand{\prf}{ \begin{proof} }
\newcommand{\epr}{ \end{proof} }

\newcommand{\nc}{\newcommand}
\newcommand{\iso}{\stackrel{\sim}{\longrightarrow}}
\nc{\HC}{{\mathcal{HC}}}
\nc{\on}{\operatorname}
\nc{\BA}{{\mathbb{A}}}
\nc{\BC}{{\mathbb{C}}}
\nc{\BG}{{\mathbb{G}}}
\nc{\BM}{{\mathbb{M}}}
\nc{\BN}{{\mathbb{N}}}
\nc{\BQ}{{\mathbb{Q}}}
\nc{\BP}{{\mathbb{P}}}
\nc{\BR}{{\mathbb{R}}}
\nc{\BZ}{{\mathbb{Z}}}
\nc{\BS}{{\mathbb{S}}}

\nc{\CA}{{\mathcal{A}}}
\nc{\CB}{{\mathcal{B}}}
\nc{\CalC}{{\mathcal C}}
\nc{\CalD}{{\mathcal D}}
\nc{\CE}{{\mathcal{E}}}
\nc{\CF}{{\mathcal{F}}}
\nc{\CG}{{\mathcal{G}}}
\nc{\CH}{{\mathcal{H}}}
\nc{\CK}{{\mathcal{K}}}
\nc{\CL}{{\mathcal{L}}}
\nc{\CM}{{\mathcal{M}}}
\nc{\CMM}{{\mathcal{M}^{\operatorname{gen}}_\hbar(-\rho)}}
\nc{\CN}{{\mathcal{N}}}
\nc{\CO}{{\mathcal{O}}}
\nc{\CP}{{\mathcal{P}}}
\nc{\CQ}{{\mathcal{Q}}}
\nc{\CR}{{\mathcal{R}}}
\nc{\CS}{{\mathcal{S}}}
\nc{\CT}{{\mathcal{T}}}
\nc{\CU}{{\mathcal{U}}}
\nc{\CV}{{\mathcal{V}}}
\nc{\CW}{{\mathcal{W}}}
\nc{\CX}{{\mathcal{X}}}
\nc{\CY}{{\mathcal{Y}}}
\nc{\CZ}{{\mathcal{Z}}}

\nc{\gen}{{\operatorname{gen}}}
\nc{\cM}{{\check{\mathcal M}}{}}
\nc{\csM}{{\check{\mathcal A}}{}}
\nc{\obM}{{\overset{\circ}{\mathbf M}}{}}
\nc{\oCA}{{\overset{\circ}{\mathcal A}}{}}
\nc{\obA}{{\overset{\circ}{\mathbf A}}{}}
\nc{\ooM}{{\overset{\circ}{M}}{}}
\nc{\osM}{{\overset{\circ}{\mathsf M}}{}}
\nc{\vM}{{\overset{\bullet}{\mathcal M}}{}}
\nc{\nM}{{\underset{\bullet}{\mathcal M}}{}}
\nc{\obD}{{\overset{\circ}{\mathbf D}}{}}
\nc{\cp}{{\overset{\circ}{\mathbf p}}{}}
\nc{\ofZ}{{\overset{\circ}{\mathfrak Z}}{}}
\nc{\oCZ}{{\overset{\circ}{\mathcal Z}}{}}

\nc{\fa}{{\mathfrak{a}}}
\nc{\fb}{{\mathfrak{b}}}
\nc{\fc}{{\mathfrak{c}}}
\nc{\fd}{{\mathfrak{d}}}
\nc{\ff}{{\mathfrak{f}}}
\nc{\fg}{{\mathfrak{g}}}
\nc{\fgl}{{\mathfrak{gl}}}
\nc{\fh}{{\mathfrak{h}}}
\nc{\fj}{{\mathfrak{j}}}
\nc{\fm}{{\mathfrak{m}}}
\nc{\fn}{{\mathfrak{n}}}
\nc{\fu}{{\mathfrak{u}}}
\nc{\fp}{{\mathfrak{p}}}
\nc{\frr}{{\mathfrak{r}}}
\nc{\fs}{{\mathfrak{s}}}
\nc{\ft}{{\mathfrak{t}}}
\nc{\fy}{{\mathfrak{y}}}
\nc{\ofT}{{\overline{\mathfrak T}}}
\nc{\ofS}{{\overline{\mathfrak S}}}
\nc{\fsl}{{\mathfrak{sl}}}
\nc{\hsl}{{\widehat{\mathfrak{sl}}}}
\nc{\hgl}{{\widehat{\mathfrak{gl}}}}
\nc{\hg}{{\widehat{\mathfrak{g}}}}
\nc{\chg}{{\widehat{\mathfrak{g}}}{}^\vee}
\nc{\hn}{{\widehat{\mathfrak{n}}}}
\nc{\chn}{{\widehat{\mathfrak{n}}}{}^\vee}

\nc{\fA}{{\mathfrak{A}}}
\nc{\fB}{{\mathfrak{B}}}
\nc{\fC}{{\mathfrak{C}}}
\nc{\fD}{{\mathfrak{D}}}
\nc{\fE}{{\mathfrak{E}}}
\nc{\fF}{{\mathfrak{F}}}
\nc{\fG}{{\mathfrak{G}}}
\nc{\fI}{{\mathfrak{I}}}
\nc{\fJ}{{\mathfrak{J}}}
\nc{\fK}{{\mathfrak{K}}}
\nc{\fL}{{\mathfrak{L}}}
\nc{\fM}{{\mathfrak{M}}}
\nc{\fN}{{\mathfrak{N}}}
\nc{\frP}{{\mathfrak{P}}}
\nc{\fQ}{{\mathfrak Q}}
\nc{\fR}{{\mathfrak R}}
\nc{\fS}{{\mathfrak S}}
\nc{\fT}{{\mathfrak{T}}}
\nc{\fU}{{\mathfrak{U}}}
\nc{\fV}{{\mathfrak{V}}}
\nc{\fW}{{\mathfrak{W}}}
\nc{\fZ}{{\mathfrak{Z}}}

\nc{\bb}{{\mathbf{b}}}
\nc{\bc}{{\mathbf{c}}}
\nc{\be}{{\mathbf{e}}}
\nc{\bj}{{\mathbf{j}}}
\nc{\bn}{{\mathbf{n}}}
\nc{\bp}{{\mathbf{p}}}
\nc{\bq}{{\mathbf{q}}}
\nc{\bs}{{\mathbf{s}}}
\nc{\bt}{{\mathbf{t}}}
\nc{\bv}{{\mathbf{v}}}
\nc{\bx}{{\mathbf{x}}}
\nc{\by}{{\mathbf{y}}}
\nc{\bw}{{\mathbf{w}}}
\nc{\bA}{{\mathbf{A}}}
\nc{\bB}{{\mathbf{B}}}
\nc{\bC}{{\mathbf{C}}}
\nc{\bK}{{\mathbf{K}}}
\nc{\bD}{{\mathbf{D}}}
\nc{\bH}{{\mathbf{H}}}
\nc{\bI}{{\mathbf{I}}}
\nc{\bM}{{\mathbf{M}}}
\nc{\bN}{{\mathbf{N}}}
\nc{\bO}{{\mathbf{O}}}
\nc{\bQ}{{\mathbf Q}}
\nc{\bS}{{\mathbf{S}}}
\nc{\bT}{{\mathbf{T}}}
\nc{\bU}{{\mathbf{U}}}
\nc{\bV}{{\mathbf{V}}}
\nc{\bW}{{\mathbf{W}}}
\nc{\bX}{{\mathbf{X}}}
\nc{\bP}{{\mathbf{P}}}
\nc{\bZ}{{\mathbf{Z}}}

\nc{\sa}{{\mathsf{a}}}
\nc{\sfb}{{\mathsf{b}}}
\nc{\sA}{{\mathsf{A}}}
\nc{\sB}{{\mathsf{B}}}
\nc{\sC}{{\mathsf{C}}}
\nc{\sD}{{\mathsf{D}}}
\nc{\sF}{{\mathsf{F}}}
\nc{\sK}{{\mathsf{K}}}
\nc{\sM}{{\mathsf{M}}}
\nc{\sO}{{\mathsf{O}}}
\nc{\sQ}{{\mathsf{Q}}}
\nc{\sP}{{\mathsf{P}}}
\nc{\sR}{{\mathsf{R}}}
\nc{\sS}{{\mathsf{S}}}
\nc{\sT}{{\mathsf{T}}}
\nc{\sV}{{\mathsf{V}}}
\nc{\sW}{{\mathsf{W}}}
\nc{\sX}{{\mathsf{X}}}
\nc{\sZ}{{\mathsf{Z}}}
\nc{\sfp}{{\mathsf{p}}}
\nc{\sr}{{\mathsf{r}}}
\nc{\st}{{\mathsf{t}}}
\nc{\sv}{{\mathsf{v}}}
\nc{\sy}{{\mathsf{y}}}
\nc{\sfc}{{\mathsf{c}}}
\nc{\sd}{{\mathsf{d}}}
\nc{\sg}{{\mathsf{g}}}
\nc{\sfl}{{\mathsf{l}}}

\nc{\BK}{{\bar{K}}}

\nc{\tA}{{\widetilde{\mathbf{A}}}}
\nc{\tB}{{\widetilde{\mathcal{B}}}}
\nc{\tg}{{\widetilde{\mathfrak{g}}}}
\nc{\tG}{{\widetilde{G}}}
\nc{\TM}{{\widetilde{\mathbb{M}}}{}}
\nc{\tO}{{\widetilde{\mathsf{O}}}{}}
\nc{\tU}{{\widetilde{\mathfrak{U}}}{}}
\nc{\TZ}{{\tilde{Z}}}
\nc{\tZ}{\widetilde{Z}{}}
\nc{\tx}{{\tilde{x}}}
\nc{\tbv}{{\tilde{\bv}}}
\nc{\tfP}{{\widetilde{\mathfrak{P}}}{}}
\nc{\tz}{{\tilde{\zeta}}}
\nc{\tmu}{{\tilde{\mu}}}

\nc{\td}{\ddot{\underline{d}}{}}
\nc{\tzeta}{\widetilde{\zeta}{}}
\nc{\hd}{{\widehat{\underline{d}}}}
\nc{\hG}{{\widehat{G}}}
\nc{\hBP}{\widehat{\mathbb P}{}}
\nc{\hQ}{{\widehat{Q}}}
\nc{\hsM}{\widehat{\mathsf M}{}}
\nc{\hfM}{\widehat{\mathfrak M}{}}
\nc{\hCP}{\widehat{\mathcal P}{}}
\nc{\hCR}{\widehat{\mathcal R}{}}
\nc{\hCS}{{\widehat{\mathcal S}}}
\nc{\hfZ}{\widehat{\mathfrak Z}{}}

\nc{\urho}{\underline{\rho}}
\nc{\uB}{\underline{B}}
\nc{\uC}{{\underline{\mathbb{C}}}}
\nc{\ui}{\underline{i}}
\nc{\ofP}{{\overline{\mathfrak{P}}}}
\nc{\oA}{\overset{\circ}{\mathbb A}{}}

\nc{\hrho}{{\hat{\rho}}}

\nc{\unl}{\underline}
\nc{\ol}{\overline}
\nc{\one}{{\mathbf{1}}}
\nc{\two}{{\mathbf{t}}}

\nc{\Tot}{{\mathop{\operatorname{\rm Tot}}}}
\nc{\Hilb}{{\mathop{\operatorname{\rm Hilb}}}}
\nc{\End}{{\operatorname{End}}}
\nc{\Ext}{{\operatorname{Ext}}}
\nc{\Hom}{{\operatorname{Hom}}}
\nc{\Sym}{{\operatorname{Sym}}}
\nc{\CHom}{{\mathop{\operatorname{{\mathcal{H}}om}}}}
\nc{\defi}{{\mathop{\operatorname{\rm def}}}}
\nc{\length}{{\mathop{\operatorname{\rm length}}}}
\nc{\Bun}{{\on{Bun}}}
\nc{\Cliff}{{\mathsf{Cliff}}}
\nc{\bGr}{{\mathbf{Gr}}}
\nc{\bFl}{{\mathbf{Fl}}}
\nc{\bIw}{{\mathbf{Iw}}}
\nc{\Fib}{{\mathsf{Fib}}}
\nc{\Coh}{{\mathsf{Coh}}}
\nc{\FCoh}{{\mathsf{FCoh}}}

\nc{\reg}{{\on{reg}}}
\nc{\rat}{{\on{rat}}}
\nc{\trig}{{\on{trig}}}

\nc{\cplus}{{\mathbf{C}_+}}
\nc{\cminus}{{\mathbf{C}_-}}
\nc{\cthree}{{\mathbf{C}_*}}
\nc{\Qbar}{{\bar{Q}}}
\nc{\Fl}{{{\mathcal F}\ell}}
\nc{\Gr}{{\on{Gr}}}
\nc{\Iw}{{\on{Iw}}}
\nc{\bh}{{\bar{h}}}
\nc{\bOmega}{{\overline{\Omega}}}
\nc\tGr{\widetilde{\Gr}}
\nc{\ul}{\underline}
\nc{\seq}[1]{\stackrel{#1}{\sim}}
\nc\ogu{\overline{G/U}}
\nc\chlam{\check{\lam}}
\nc\St{\operatorname{St}}
\nc{\bLambda}{{\boldsymbol{\Lambda}}}
\nc\uS{\underline{S}}
\nc\QM{\mathcal{QM}}
\nc{\hQM}{\widehat{\mathcal{QM}}{}}
\nc{\dQM}{{}^{t_0}\!{\mathcal{QM}}}
\nc{\chmu}{\check{\mu}}
\nc{\CHH}{{\CH\!\!\CH}}

\nc{\hZ}{\,^\dagger\!Z}
\nc{\dZ}{\,^\dagger\!\overset{\circ}{Z}{}}
\nc{\buZ}{\,^\dagger\!\overset{\bullet}{Z}{}}
\nc{\oY}{\overset{\circ}{Y}{}}
\nc{\oZ}{\overset{\circ}{Z}{}}
\nc{\Cx}{{C^\dagger}}
\nc{\pr}{{\on{pr}}}
\nc{\gaff}{{\fg_{\on{aff}}}}

\nc{\Lam}{\Lambda}
\nc{\calB}{\mathcal B}
\nc{\CC}{\mathbb C}
\nc{\ZZ}{\mathbb Z}
\nc{\alp}{\alpha}
\nc{\PP}{\mathbb P}
\nc\x{\times}
\renewcommand{\o}{{\check\omega}}
\nc\lam{\lambda}
\nc\gam{\gamma}
\nc\calF{\mathcal F}
\nc\disj{\operatorname{disj}}

\newcommand{\slhat}{\widehat {\mathfrak s \mathfrak l}}
\renewcommand{\det}{{\mathrm d \mathrm e \mathrm t} \ }

\usepackage{amsmath}
\usepackage{kbordermatrix}

\begin{document}

\bigskip

\author{%Three men in a boat
Michael Finkelberg, Alexander Kuznetsov and Leonid Rybnikov \\(with an appendix by Galyna Dobrovolska)}
\title{Towards a cluster structure on trigonometric zastava}

\dedicatory{To Sasha Beilinson with love and gratitude}

\begin{abstract}
%We clusterize.
We study a moduli problem on a nodal curve of arithmetic genus 1, whose solution
is an open subscheme in the zastava space for projective line. This moduli space
is equipped with a natural Poisson structure, and we compute it in a natural
coordinate system. We compare this Poisson structure with the trigonometric
Poisson structure on the transversal slices in an affine flag variety.
We conjecture that certain generalized minors give rise to a cluster structure
on the trigonometric zastava.
\end{abstract}

\maketitle %\footnote{Note to the author: everyone has his own psychological
%problems,  some at least make an effort not to put them on display.}
%\tableofcontents

\sec{Intro}{Introduction}

\ssec{rati}{Zastava and euclidean monopoles}
Let $G$ be an almost simple simply connected algebraic group over $\CC$. We denote by $\calB$ the flag variety of
$G$. Let us also fix a pair of opposite Borel subgroups $B$, $B_-$ whose intersection is a maximal torus $T$. Let $\Lam$ denote the cocharacter lattice of $T$; since $G$ is assumed to be simply connected, this is also the coroot lattice of $G$.
We denote by $\Lam_+\subset \Lam$ the sub-semigroup spanned by positive coroots.

It is well-known that $H_2(\calB,\ZZ)=\Lam$ and that an element $\alp\in H_2(\calB,\ZZ)$ is representable by an effective algebraic curve
if and only if $\alp\in \Lam_+$.
Let $\oZ^{\alp}$ denote the space of maps $C=\PP^1\to \calB$ of degree $\alp$ sending $\infty\in \PP^1$ to $B_-\in \calB$.
It is known~\cite{fkmm} that this is a smooth symplectic affine algebraic variety, which can be identified with the hyperk\"ahler moduli space of
framed $G$-monopoles on $\BR^3$ with maximal symmetry breaking at infinity of charge $\alp$~\cite{j},~\cite{j'}.

The monopole space $\oZ^{\alp}$ has a natural partial compactification
$Z^{\alp}$ ({\em zastava} scheme). It can be realized as the moduli
space of based {\em quasi-maps} of degree $\alp$; set-theoretically it can be described in the following way:
$$
Z^{\alp}=\bigsqcup\limits_{0\leq\beta\leq \alp} \oZ^{\beta}\x
\BA^{\alpha-\beta},
$$
where for $\gam\in \Lam_+$ we denote by $\BA^{\gam}$ the space of all colored divisors $\sum\gam_i x_i$ with
$x_i\in \BA^1$, $\gam_i\in \Lam_+$ such that $\sum \gam_i=\gam$.

The zastava space is equipped with a {\em factorization} morphism
$\pi_\alpha\colon Z^\alpha\to\BA^\alpha$ whose restriction to
$\oZ^\alpha\subset Z^\alpha$ has a simple geometric meaning: for a based map
$\varphi\in\oZ^\alpha$ the colored divisor $\pi_\alpha(\phi)$ is just the
pullback of the colored Schubert divisor $D\subset\CB$ equal to the complement
of the open $B$-orbit in $\CB$. The morphism
$\pi_\alpha\colon \oZ^\alpha\to\BA^\alpha$ is the {\em Atiyah-Hitchin} integrable
system (with respect to the above symplectic structure): all the fibers of
$\pi_\alpha$ are Lagrangian.

A system of \'etale birational coordinates on $\oZ^{\alp}$ was introduced
in~\cite{fkmm}. Let us recall the definition
for $G=SL(2)$. In this case $\oZ^{\alp}$ consists of all maps $\PP^1\to \PP^1$ of degree $\alp$ which send
$\infty$ to $0$. We can represent such a map by a rational function $\frac{R}{Q}$ where $Q$ is a monic polynomial of degree
$\alp$ and $R$ is a polynomial of degree $<\alp$. Let $w_1,\ldots,w_{\alp}$ be the zeros of $Q$. Set $y_r=R(w_r)$.
Then the functions $(y_1,\ldots,y_\alp,w_1,\ldots, w_{\alp})$ form a system of \'etale birational coordinates on $\oZ^{\alp}$, and the above mentioned
symplectic form in these coordinates reads
$\Omega_\rat=\sum_{r=1}^\alpha\frac{dy_r\wedge dw_r}{y_r}$.

For general $G$ the definition of the above coordinates is quite similar. In this case given a point in $\oZ^{\alp}$ we can define polynomials $R_i,Q_i$ where $i$
runs through the set $I$ of vertices of the Dynkin diagram of
$G,\ \alpha=\sum a_i\alpha_i$, and

(1) $Q_i$ is a monic polynomial of degree $a_i$,

(2) $R_i$ is a polynomial of degree $<a_i$.

Hence, we can define (\'etale, birational) coordinates
$(y_{i,r},w_{i,r})$ where $i\in I$ and
$r=1,\ldots,a_i$. Namely, $w_{i,r}$ are the roots of $Q_i$, and
$y_{i,r}=R_i(w_{i,r})$. The Poisson brackets of these coordinates
with respect to the above symplectic form are as follows:
$\{w_{i,r},w_{j,s}\}_\rat=0,\
\{w_{i,r},y_{j,s}\}_\rat=\check{d}_i\delta_{ij}\delta_{rs}y_{j,s},\
\{y_{i,r},y_{j,s}\}_\rat=
(\check\alpha_i,\check\alpha_j)\frac{y_{i,r}y_{j,s}}{w_{i,r}-w_{j,s}}$ for
$i\ne j$, and finally $\{y_{i,r},y_{i,s}\}_\rat=0$. Here $\check\alpha_i$ is a
simple root, $(,)$ is the invariant scalar product on $(\on{Lie}T)^*$ such
that the square length of a short root is 2, and
$\check{d}_i=(\check\alpha_i,\check\alpha_i)/2$.

Now recall that the standard {\em rational} $r$-matrix for $\fg=\on{Lie}G$
gives rise to a Lie bialgebra structure on $\fg[z^{\pm1}]$ corresponding to
the Manin triple $\fg[z],z^{-1}\fg[z^{-1}],\fg[z^{\pm1}]$. This in turn
gives rise to a Poisson structure on the affine Grassmannian
$\Gr_G=G[z^{\pm1}]/G[z]$. The transversal slices $\CW^\lambda_\mu$ from
a $G[z]$-orbit $\Gr_G^\mu$ to another orbit $\Gr_G^\lambda$ (here
$\mu\leq\lambda$ are dominant coweights of $G$) are examples of symplectic
leaves of the above Poisson structure. According to~\cite{bf},
the zastava spaces are ``stable limits" of the above slices. More precisely,
for $\alpha=\lambda-\mu$ there is a birational Poisson map
$s^{\lambda^*}_{\mu^*}\colon \CW^{\lambda^*}_{\mu^*}\to Z^\alpha$ (here
$\lambda^*:=-w_0\lambda,\ \mu^*=-w_0\mu$, and $w_0$ is the longest element
of the Weyl group $W=W(G,T)$).

\ssec{tr}{Trigonometric zastava and periodic monopoles}
We have an open subset $\BG_m^\alpha\subset\BA^\alpha$ (colored divisors not
meeting $0\in\BA^1$), and we introduce the open subscheme of
{\em trigonometric zastava}
$\hZ^\alpha:=\pi_\alpha^{-1}\BG_m^\alpha\subset Z^\alpha$, and its smooth
open affine subvariety of {\em periodic monopoles}
$\dZ^\alpha:=\hZ^\alpha\cap\oZ^\alpha$. These schemes are solutions of the
following modular problems.

Let $\Cx$ be an irreducible nodal curve of arithmetic genus 1 obtained by
gluing the points $0,\infty\in C=\BP^1$, so that $\pi\colon C\to\Cx$ is the
normalization. Let $\fc\in\Cx$ be the singular point.
The moduli space $\Bun_T^0(\Cx)$ of $T$-bundles on $\Cx$
of degree 0 is canonically identified with the Cartan torus $T$ itself.
We fix a $T$-bundle $\CF_T$
which corresponds to a {\em regular} point $t\in T^\reg$. Then $\dZ^\alpha$
is the moduli space of the following data:\\ (a) a trivialization $\tau_\fc$
of the fiber of $\CF_T$ at the singular point $\fc\in\Cx$;\\
(b) a $B$-structure $\phi$ in the induced $G$-bundle
$\CF_G=\CF_T\overset{T}{\times}G$ of degree $\alpha$
which is transversal to $\CF_B=\CF_T\overset{T}{\times}B$ at $\fc$.\\
The scheme $\hZ^\alpha$ is the moduli space of the similar data with the only
difference: we allow a $B$-structure in (b) to be {\em generalized} (i.e.
to acquire defects at certain points of $\Cx$), but
require it to have no defect at $\fc$.

As a regular Cartan element $t$ varies, the above moduli spaces become fibers
of a single family.
More precisely, we consider the following moduli problem:\\
(t) a $T$-bundle $\CF_T$ of degree 0 on $\Cx$ corresponding to a regular
element of $T$;\\ (a,b) as above;\\ (c) a trivialization $f_\fc$ at $\fc$ of the $T$-bundle $\phi_T$ induced from the $B$-bundle $\phi$ in (b).\\
This moduli problem is represented by a scheme $\oY^\alpha\subset Y^\alpha$
(depending on whether the $B$-structure in (b) is genuine or generalized).
Note that $Y^\alpha$ is equipped with an action of $T\times T$ changing the
trivializations in (a,c). We prove that $\oY^\alpha$ is a smooth affine
variety equipped with a natural projection
$\varpi\colon \oY^\alpha\to\dZ^\alpha$, and we construct a nondegenerate bivector
field on $\oY^\alpha$ arising from a differential in a spectral
sequence involving the tangent and cotangent bundles of $\oY^\alpha$
(this construction is a trigonometric degeneration of the
construction~\cite{fo} for elliptic curves; its rational analogue was worked
out in~\cite{fkmm}).
This bivector field descends to $\dZ^\alpha$ under the projection
$\varpi\colon \oY^\alpha\to\dZ^\alpha$ and gives rise to a nondegenerate Poisson
structure, i.e.\ a symplectic form on $\dZ^\alpha$.
%The projection $\oY^\alpha\to T^\reg$ plays a role of the moment map for one of
%the above $T$-actions, and $\dZ^\alpha$ is nothing but the quasihamiltonian
%reduction of $\oY^\alpha$.
It would be interesting to obtain this symplectic form by the method 
of~\cite[4.2]{s}.
The Poisson brackets of the coordinates on $\dZ^\alpha$ are as follows:
$\{w_{i,r},w_{j,s}\}_\trig=0,\
\{w_{i,r},y_{j,s}\}_\trig=\check{d}_i\delta_{ij}\delta_{rs}w_{j,s}y_{j,s},\ \{y_{i,r},y_{j,s}\}_\trig=
(\check\alpha_i,\check\alpha_j)\frac{(w_{i,r}+w_{j,s})y_{i,r}y_{j,s}}
{2(w_{i,r}-w_{j,s})}$ for
$i\ne j$, and finally $\{y_{i,r},y_{i,s}\}_\trig=0$.
In particular, the projection $\pi_\alpha\colon \dZ^\alpha\to\BG_m^\alpha$ is the
{\em trigonometric Atiyah-Hitchin} integrable system (for $G=SL(2)$ this
system goes back at least to~\cite{fg}).

Now recall that the standard {\em trigonometric} $r$-matrix for $\fg$ gives
rise to a Lie bialgebra structure on $\fg((z^{-1}))\oplus\ft$ which in turn gives rise
to a Poisson structure on the affine flag variety $\Fl_G$ (the quotient of
$G[z^{\pm1}]$ with respect to an Iwahori subgroup).
%The symplectic leaves of this Poisson structure are
The intersections $\Fl^w_y$ of the opposite
Iwahori orbits (aka {\em open Richardson varieties}) are Poisson subvarieties
of $\Fl_G$. Here $w,y$ are elements
of the affine Weyl group $W_a=W\ltimes\Lambda$. For dominant coweights
$\mu\leq\lambda\in\Lambda$ such that $\lambda-\mu=\alpha$, 
and the longest element $w_0$ of the finite Weyl
group $W$, the Richardson variety $\Fl^{w_0\times\lambda^*}_{w_0\times\mu^*}$
is a symplectic leaf of $\Fl_G$, the projection
$\on{pr}\colon \Fl^{w_0\times\lambda^*}_{w_0\times\mu^*}\to\CW^{\lambda^*}_{\mu^*}$ is an
open embedding, and the composition
$s^{\lambda^*}_{\mu^*}\circ\on{pr}\colon \Fl^{w_0\times\lambda^*}_{w_0\times\mu^*}\to
Z^\alpha$ is a {\em symplectomorphism} onto its image
$\dZ^\alpha\subset Z^\alpha$ equipped with the trigonometric symplectic
structure.

\ssec{clu}{Cluster aspirations}
It seems likely that the construction due to~B.~Leclerc~\cite{l} extends from
the open Richardson varieties in the type ADE finite flag varieties to
the case of the affine flag varieties, and provides $\Fl^w_y$ with a cluster
structure (even in the nonsimply laced case, cf.~\cite{ley}).
This structure can be transferred from $\Fl^{w_0\times\lambda^*}_{w_0\times\mu^*}$
to $\dZ^\alpha$ via the above symplectomorphism.
If $\alpha$ is a dominant coweight of $G$, a reduced decomposition of
$w_0\times\mu^*$ is the beginning of a reduced decomposition of
$w_0\times\lambda^*$, and the existence of cluster structure on
$\Fl^{w_0\times\lambda^*}_{w_0\times\mu^*}$ is known for arbitrary symmetric affine
Kac-Moody algebra.
In case of $G=SL(2)$ the resulting cluster structure
on the moduli space of periodic monopoles was discovered in~\cite{agsv},
which served as the starting point of the present note. It seems likely that
for general $G$ the Gaiotto-Witten superpotential on $Z^\alpha$
(see e.g.~\cite{bdf}) restricted to $\dZ^\alpha$ is totally positive in the
above cluster structure (see~\refs{clus} for details).

\rem{lecler}
Implicit in the above discussion when $\alpha$ is {\em dominant} (as a coweight
of $G$) is an affine open embedding
$\dZ^\alpha\subset\oZ^\alpha\hookrightarrow\BA^{2|\alpha|}$ into an affine space. 
Indeed, in this case the lengths of the affine Weyl group elements satisfy
$\ell(w_0\times\mu^*)+\ell(\alpha^*)=\ell(w_0\times\lambda^*)$, and we are
in the situation of~\cite[Section~5]{l}; hence according 
to {\em loc.~cit.}~and~\cite[Theorem~2.12]{l}, 
$\Fl^{w_0\times\lambda^*}_{w_0\times\mu^*}$ is an open
subvariety of an affine space.
Here is a modular interpretation of the above open embedding: 
$\BA^{2|\alpha|}$ is the moduli space
of $B$-bundles $\phi_B$ on $\BP^1$ equipped with
a trivialization $(\phi_B)_\infty\iso B$ of the fiber at $\infty\in\BP^1$,
such that the induced $T$-bundle
(under projection $B\twoheadrightarrow T$) has degree $\alpha$.
\erem

\ssec{coulomb}{Relation to Coulomb branches of $4d\ \CN=2$ quiver gauge
theories}
Let $G=SL(2)$, and $\alpha=a\in\BN$. Then the methods of~\cite{bfm}
establish an isomorphism $\BC[\dZ^a]\simeq K^{GL(a)_\CO}(\Gr_{GL(a)})$
(where $\CO=\BC[[z]]$). More generally, for $G$ of type $ADE$ let us orient
its Dynkin diagram, and for $\alpha=\sum_{i\in I}a_i\alpha_i$ let us consider
a representation $V$ of the Dynkin quiver such that $\dim V_i=a_i$.
The group $GL(V):=\prod_{i\in I}GL(V_i)$ acts in
$\bN=\bigoplus_{i\to j}\Hom(V_i,V_j)$. Following~\cite{BFN1} we consider the
moduli space $\CR$ of triples $(\CF,\sigma,s)$ where $\CF$ is a $GL(V)$-bundle
on the formal disc $D=\on{Spec}\BC[[z]],\ \sigma$ is its trivialization over
the punctured disc $D^*=\on{Spec}\BC((z))$, and $s$ is a section of the
associated vector bundle $\CF\times^G\bN$ such that it is sent to a regular
section of the trivial bundle under $\sigma$. The group $GL(V)_\CO$ acts
naturally on $\CR$, and as in~\cite{BFN1} one can define the equivariant
$K$-theory $K^{GL(V)_\CO}(\CR)$ and equip it with the convolution algebra
structure. Moreover, as in~\cite{BFN3} one can establish an isomorphism
$\BC[\dZ^\alpha]\simeq K^{GL(V)_\CO}(\CR)$ such that the factorization morphism
$\pi_\alpha\colon \dZ^\alpha\to\BG_m^\alpha$ corresponds to the embedding of the
equivariant $K$-theory of the point:
$K^{GL(V)_\CO}(pt)\hookrightarrow K^{GL(V)_\CO}(\CR)$.

Yet more generally, given a framing $W_i,\ \dim W_i=l_i$, we set $\lambda=\sum l_i\omega_i$
and consider a representation $\bN'=\bN\oplus\bigoplus_i\Hom(W_i,V_i)$ of $GL(V)$. It gives
rise to the space $\CR'$ of triples as above, and one can prove as in~\cite{BFN3} that
the convolution algebra $K^{GL(V)_\CO}(\CR')$ is isomorphic to the coordinate ring
of the moduli space $^\dagger\ol\CW{}^{\lambda^*}_{\mu^*}$ of the triples
$(\CP,\sigma,\phi)$ where $\CP$ is a $G$-bundle on $C;\ \sigma$ is a trivialization
of $\CP$ off $1\in C$ having a pole of degree $\leq\lambda^*$ at $1\in C$, and
$\phi$ is a $B$-structure on $\CP$ of degree $-\mu$ having the fiber $B_-$ at
$\infty\in C$ and transversal to $B$ at $0\in C$ (a trigonometric slice).
Note that the relation of $^\dagger\ol\CW{}^{\lambda^*}_{\mu^*}$ to the 
Richardson varieties $\Fl^w_y$ of~\refss{tr} is unclear since the former
``knows'' about 3 points $0,1,\infty\in\BP^1$ while the latter only ``knows''
about 2 points.

\ssec{contents}{Contents}
In~\refs{node} the moduli spaces $\dZ^\alpha\subset\hZ^\alpha$ and
$\oY^\alpha\subset Y^\alpha$ are introduced. In~\refs{sympl} a bivector field
on $\oY^\alpha$ is introduced and the resulting Poisson bracket $\{,\}_\trig$
on $\dZ^\alpha$ is computed in coordinates $(w_{i,r},y_{i,r})$. The
technicalities of this computation occupy the bulk of the present note.
In~\refs{slice} we compare the rational Poisson bracket on the slices
$\CW^{\lambda^*}_{\mu^*}\subset\Gr_G$ with the Poisson bracket $\{,\}_\rat$ on the
monopole moduli space $\oZ^\alpha$. We also compare the trigonometric Poisson
bracket on the slices $\Fl^{w_0\times\lambda^*}_{w_0\times\mu^*}\subset\Fl_G$
with the Poisson bracket $\{,\}_\trig$ on the periodic monopole moduli
space $\dZ^\alpha$. \refs{clus} contains a few comments on the
cluster structures. Finally, the Appendix~\refs{append} by Galyna Dobrovolska
identifies our cluster structure on $\dZ^\alpha$ for $G=SL_2$ with the one
of~\cite{agsv}.

\ssec{ack}{Acknowledgments}
We are grateful to the organizers of the 10th Lunts Summer School ``Cluster
Algebras" where this note was conceived. We are much obliged to
R.~Bezrukavnikov, A.~Braverman, B.~Feigin, V.~Fock, D.~Gaiotto, M.~Gekhtman, 
J.~Kamnitzer, B.~Leclerc, H.~Nakajima, M.~Shapiro and M.~Yakimov for
their generous and patient explanations.
%The financial support from the Government of the Russian Federation within the
%framework of the implementation of the 5-100 Programme Roadmap of the National
%Research University  Higher School of Economics, AG Laboratory  is acknowledged
%by M.F and A.K. The financial support from the Government of the Russian
%Federation within the framework of the implementation of the 5-100 Programme
%Roadmap of the National Research University  Higher School of Economics, RTMP
%Laboratory  is acknowledged by L.R.
This work has been funded by the Russian Academic Excellence
Project `5-100'.
L.R. was supported by Russian President grant MK-2121.2014.1.

\sec{node}{Trigonometric zastava}

\ssec{Cdag}{Line bundles on a nodal curve $C^\dagger$}
Let $\Cx$ be an irreducible nodal curve of arithmetic genus 1, and let
$\pi\colon C\to\Cx$ be its normalization. Then $C$ is a projective line.
We fix a coordinate function $z$ on $C$ such that the preimage of the node
$\fc\in\Cx$ consists of the points $0,\infty\in C$. We have
$\on{Pic}^0(\Cx)=\BG_m$: any line bundle on $\Cx$ is obtained by descent from
the one on $C$ gluing its fibers at $0$ and $\infty$; a degree 0 line bundle
on $C$ is trivial, so its fibers at $0$ and $\infty$ are canonically identified.
Moreover, the above choice of a coordinate $z$ on $C$ gives rise to an
identification $\on{Pic}^n(\Cx)\cong\BG_m$ for any $n\in\BZ$: if $n>0$, and
$s$ is a section of $\CL\in\on{Pic}^n(\Cx)$ not vanishing at $\fc$, then
$\on{div}(s)\in\on{Sym}^n(\Cx\setminus\fc)=\on{Sym}^n(\BG_m)$ (this
identification makes use of the coordinate $z$); we have a multiplication
morphism $m\colon \on{Sym}^n(\BG_m)\to\BG_m$, and finally $\CL\mapsto
m(\on{div}(s))\in\BG_m$ (the result is independent of the choice of a section
$s$). If $\CL\in\on{Pic}^n(\Cx)$ for $n<0$, then $\CL$ goes to the inverse
of the class of $\CL^{-1}$.

For the canonical line bundle $\omega_\Cx$, we have the following exact
sequence:
\eq{rho_+}
0\to\omega_\Cx\to\pi_*\omega_C(\{0\}+\{\infty\})
\stackrel{\rho_+}{\longrightarrow}\BC_\fc\to0
\end{equation}
where
$\rho_+(\xi)=\on{Res}_0(\xi)+\on{Res}_\infty(\xi)$. The line bundle $\omega_\Cx$
is trivial, with trivializing section $z^{-1}dz$. In what follows we will
freely use the above identification $\omega_\Cx\cong\CO_\Cx$.

We define the {\em theta-characteristic} $\theta\in\on{Pic}^0(\Cx)$ as a unique
{\em nontrivial} line bundle such that $\theta^2=\omega_\Cx$. It enters the
following exact sequence:
\eq{rho_-}
0\to\theta\to\pi_*\omega_C(\{0\}+\{\infty\})
\stackrel{\rho_-}{\longrightarrow}\BC_\fc\to0
\end{equation}
where $\rho_-(\xi)=\on{Res}_0(\xi)-\on{Res}_\infty(\xi)$.

From the above two sequences we have the natural embeddings
$\pi_*\omega_C\hookrightarrow\omega_\Cx\hookrightarrow\pi_*\omega_C(\{0\}+\{\infty\})$
and
$\pi_*\omega_C\hookrightarrow\theta\hookrightarrow\pi_*\omega_C(\{0\}+\{\infty\})$.
Noting that $\omega_C(\{0\}+\{\infty\})=\CO_C$ we combine the above embeddings into
the following exact sequence:
\eq{theta omega}
0\to\pi_*\omega_C\to\theta\oplus\omega_\Cx\to\pi_*\CO_C\to0
\end{equation}

We also have natural embeddings $\pi_*\omega_C=\pi_*\CO_C(-\{0\}-\{\infty\})
\hookrightarrow\pi_*\CO_C(-\{0\})\hookrightarrow\pi_*\CO_C$ and
$\pi_*\omega_C=\pi_*\CO_C(-\{0\}-\{\infty\})
\hookrightarrow\pi_*\CO_C(-\{\infty\})\hookrightarrow\pi_*\CO_C$.
They combine into the following exact sequence:
\eq{Xi}
0\to\pi_*\omega_C\to\pi_*\Xi\to\pi_*\CO_C\to0
\end{equation}
where $\Xi$ stands for $\CO_C(-\{0\})\oplus\CO_C(-\{\infty\})$.

\ssec{BunT}{A group $G$}
Let $G$ be an almost simple simply connected algebraic group over $\CC$. We denote by $\calB$ the flag variety of
$G$. Let us also fix a pair of opposite Borel subgroups $B$, $B_-$ whose intersection is a maximal torus $T$ (thus we have
$\calB=G/B=G/B_-$). We denote by $T^\reg\subset T$ the open subset formed
by the regular elements.

Let $\Lam$ (resp. $\Lambda^\vee$) denote the cocharacter (resp. character)
lattice of $T$; since $G$ is assumed to be simply connected, this is also the coroot lattice of $G$.
We denote by $\Lam_+\subset \Lam$ the sub-semigroup spanned by positive coroots. We say that $\alp\geq \beta$ (for $\alp,\beta\in \Lam$)
if $\alp-\beta\in\Lam_+$. The simple coroots are $\{\alpha_i\}_{i\in I}$;
the simple roots are $\{\check\alpha_i\}_{i\in I}$;
the fundamental weights are $\{\check\omega_i\}_{i\in I}$.
We consider the invariant bilinear form $(,)$ on the weight
lattice $\Lambda^\vee$ such that the square length of a {\em short} simple root
$(\check\alpha_i,\check\alpha_i)=2$. We set
$\check{d}_i:=\frac{(\check\alpha_i,\check\alpha_i)}{2}$.
We fix the Chevalley
generators $(E_i,F_i,H_i)_{i\in I}$ of $\fg$. An irreducible
$G$-module with a dominant highest weight $\check\lambda\in\Lambda^\vee_+$
is denoted $V_{\check\lambda}$; we fix its highest vector $v_{\check\lambda}$.
For a weight $\check\mu\in\Lambda^\vee$ the
$\check\mu$-weight subspace of a $G$-module $V$ is denoted $V(\check\mu)$.
Finally, $W$ is the Weyl group of $G,T$; the simple reflections are denoted
$s_i,\ i\in I$, and $w_0\in W$ is the longest element.

The identification $\on{Pic}^0(\Cx)=\BG_m$ (resp. $\on{Pic}^n(\Cx)\cong\BG_m$,
depending on the choice of coordinate $z$) of~\refss{Cdag} gives rise to the
identification
$\on{Bun}_T^0(\Cx)=T$ (resp. $\on{Bun}_T^\alpha(\Cx)\cong T$). We denote by
$\on{Bun}_T^\reg(\Cx)\subset\on{Bun}_T^0(\Cx)$ the open subset corresponding to
$T^\reg\subset T$ under the above identification.

\ssec{Yalpha}{The moduli space $Y^\alpha$}
Given a $T$-bundle $\CF_T\in\on{Bun}_T^\reg(\Cx)$ we denote by $\CF_B$ (resp.
 $\CF_G$) the corresponding induced $B$-bundle (resp. $G$-bundle).

\defe{oY}
Given $\alpha\in\Lambda_+$, we define $\oY^\alpha$ as the
 moduli space of the following data:

(a) a regular $T$-bundle $\CF_T\in\on{Bun}_T^\reg(\Cx)$;

(b) a trivialization $\tau_\fc$ of the fiber of $\CF_T$ at $\fc\in\Cx$;

(c) a $B$-structure
$\phi$ in $\CF_G$ of degree $\alpha$ (that is, the induced $T$-bundle $\phi_T$
has degree $\alpha$), such that $\phi$ is transversal to $\CF_B$ at $\fc$;

(d) a trivialization $f_\fc$ of the induced $T$-bundle $\phi_T$ at $\fc$.

We also define $Y^\alpha$ as the moduli space of the data (a--d) above where
we allow a $B$-structure in (c) to be generalized (see e.g.~\cite{bf}) but
require that it does not have a defect at $\fc\in\Cx$.
\edefe

We have a natural action of $T\times T$ on $\oY^\alpha\subset Y^\alpha$:
the first (resp. second) copy of $T$ acts via the change of trivialization
 $\tau_\fc$ (resp. $f_\fc$).

We also have a morphism
$(p,q)\colon Y^\alpha\to\on{Bun}_T^\alpha(\Cx)\times\on{Bun}_T^\reg(\Cx)\cong
T\times T^\reg$ sending $(\CF_T,\tau_\fc,\phi,f_\fc)$ to $(\phi_T,\CF_T)$.

Finally, we have a morphism $\varpi\colon Y^\alpha\to Z^\alpha$ to the zastava
space (see e.g.~\cite{bdf}) defined as follows. Recall that $Z^\alpha$ is the
moduli space of triples $(\CF^C_G,\phi_\pm^C)$ where $\CF^C_G$ is a $G$-bundle
on $C$, while $\phi_-^C$ (resp. $\phi_+^C$) is a $U$-structure in $\CF^C_G$
(resp. a generalized $B$-structure in $\CF^C_G$ of degree $\alpha$)
such that $\phi_+^C$ has no defect at $\infty\in C$, and is transversal to
$\phi_-^C$ at $\infty\in C$. Now $\varpi$ sends $(\CF_T,\tau_\fc,\phi,f_\fc)$ to
a triple $\CF^C_G:=\pi^*\CF_G,\ \phi_+^C:=\pi^*\phi$, and $\phi_-^C$ defined
as follows: $\CF^C_G$ is induced from $\CF^C_T:=\pi^*\CF_T$, and the latter
$T$-bundle is trivial and trivialized at $\infty\in\pi^{-1}(\fc)$. This 
trivialization extends uniquely to the whole of $C$, and induces a 
trivialization of $\CF^C_G$. At last, $\phi_-^C$ is a trivial $U$-structure
in the trivial $G$-bundle $\CF^C_G$ corresponding to the point $1\in G/U$.
Note that $\oY^\alpha=\varpi^{-1}(\oZ^\alpha)$ (recall that the open subset
$\oZ^\alpha\subset Z^\alpha$ is formed by the triples $(\CF^C_G,\phi_\pm^C)$
such that $\phi_+^C$ has no defect, i.e. is a usual as opposed to generalized
$B$-structure. The moduli space $\oZ^\alpha$ is isomorphic to the moduli space
of degree $\alpha$ based maps from $(C,\infty)$ to $(\CB,B_-)$).

\prop{repr} $Y^\alpha$ is represented by a scheme.
\eprop

\prf
Recall the scheme $\hQM^\alpha_\fg$ introduced in~\cite[2.3]{bf12}. It is the
moduli space of degree $\alpha$ generalized $B$-structures $\phi^C$ in the
trivial $G$-bundle on $C$, equipped with a trivialization $f_\infty$ at
$\infty\in C$ of the induced $T$-bundle $\phi^C_T$. We claim that $Y^\alpha$
is a locally closed subscheme in $T^\reg\times\hQM^\alpha_\fg$. In effect,
given a regular $T$-bundle $\CF_T\in\on{Bun}_T^\reg(\Cx)=T^\reg$, its
trivialization $\tau_\fc$ at $\fc\in\Cx$ gives rise to a trivialization
$\tau_\infty$ of $\pi^*\CF_T$ at $\infty\in C$ which extends uniquely to a
trivialization of $\pi^*\CF_T$ on $C$, and hence to a trivialization of
$\pi^*\CF_G$ on $C$. Now $\phi^C:=\pi^*\phi$ is a generalized $B$-structure
in $\pi^*\CF_G$, and the trivialization $f_\fc$ gives rise to a trivialization
$f_\infty$ of $\phi^C_T$ at $\infty\in C$. Note that $\phi^C$ has no defect
neither at $0\in C$ nor at $\infty\in C$, and its values
$\phi^C(0),\phi^C(\infty)\in\CB$ are both transversal to $B\in\CB$.
Conversely, given
$(t,\phi^C,f_\infty)\in T^\reg\times\hQM^\alpha$ such that $\phi^C$ has no
defect neither at $0\in C$ nor at $\infty\in C$, and
$\phi^C(0)=t\phi^C(\infty)\in\CB$ is transversal to $B$, we construct
$(\CF_T,\tau_\fc,\phi,f_\fc)\in Y^\alpha$ by descent from $C$ to $\Cx$.
\epr

\ssec{reduc}{A reduction of $Y^\alpha$}
Recall the {\em factorization} morphism
$\pi_\alpha\colon Z^\alpha\to\BA^\alpha=(C\setminus\{\infty\})^\alpha$
(see e.g.~\cite{bdf}). We have an open embedding
$\BG_m^\alpha=(C\setminus\{0,\infty\})^\alpha\subset\BA^\alpha$.

\defe{heart} We define the trigonometric zastava space as
$\hZ^\alpha:=\pi_\alpha^{-1}(\BG_m^\alpha)\subset Z^\alpha$.
We define the periodic monopole moduli space as
$\dZ^\alpha:=\hZ^\alpha\cap\oZ^\alpha$: a dense open smooth subscheme
of the trigonometric zastava $\hZ^\alpha$.
\edefe

Recall the action of $T\times T$ on $Y^\alpha$, and the morphism
$(p,q)\colon Y^\alpha\to T\times T^\reg$ introduced in~\refss{Yalpha}.
The action of $1\times T$ on $Y^\alpha$ is clearly free. The morphism
$\varpi\colon Y^\alpha\to Z^\alpha$ of~\refss{Yalpha} is clearly
$(1\times T)$-equivariant and gives rise to the same named morphism
$\varpi\colon Y^\alpha/(1\times T)\to Z^\alpha$.
We fix $t_0\in T^\reg$.

\prop{diamond}
For any $t_0\in T^\reg$ we have an isomorphism
$\varpi\colon q^{-1}(t_0)/(1\times T)\iso\hZ^\alpha$, and
$\varpi\colon (q^{-1}(t_0)\cap\oY^\alpha)/(1\times T)\iso\dZ^\alpha$.
\eprop

\prf
The locally closed embedding
$Y^\alpha\hookrightarrow T^\reg\times\hQM^\alpha_\fg$
constructed in the proof of~\refp{repr} gives rise to the locally closed
embedding $q^{-1}(t_0)/(1\times T)\hookrightarrow\{t_0\}\times\QM^\alpha_\fg$
where
$\QM^\alpha_\fg$ is the moduli space of degree $\alpha$ quasimaps from $C$ to
$\CB$.
More precisely, the image of $q^{-1}(t_0)/(1\times T)$ is the locally closed
subscheme $\dQM^\alpha_\fg\subset\QM^\alpha_\fg$ formed by the quasimaps that have no defects at $0,\infty\in C$,
their values $\phi^C(0)$ and $\phi^C(\infty)$ are both transversal to $B$,
and $\phi^C(0)=t_0\phi^C(\infty)$.

An open subset of $\CB$ formed by the Borels transversal to $B$ (the big
Schubert
cell) is a free orbit $U\cdot\{B_-\}$, and we will identify it with $U$ (the
unipotent radical of $B$). So for $\phi^C\in\dQM^\alpha_\fg$ we have
$\phi^C(\infty)=n\cdot\{B_-\}$ for $n\in U$, and we will simply write
$\phi^C(\infty)=n\in U$. Note that $G$ acts on $\QM^\alpha_\fg$, and $U$ acts on
$\dQM^\alpha_\fg$, and $n^{-1}\phi^C(\infty)=B_-$, i.e.
$n^{-1}\phi^C\in Z^\alpha$ is a based quasimap.

A moment of reflection shows that
$\varpi(t_0,\tau_\fc,\phi,f_\fc)=(\phi^C(\infty))^{-1}\phi^C$, and the condition
of transversality of $B$ and $\phi^C(0)$ guarantees that the value at $0\in C$
of $(\phi^C(\infty))^{-1}\phi^C$ is also transversal to $B$, i.e.
$(\phi^C(\infty))^{-1}\phi^C\in\hZ^\alpha$. Thus we have a well defined
morphism $q^{-1}(t_0)/(1\times T)\cong\dQM^\alpha_\fg\to\hZ^\alpha$, and we have
to prove that it is an isomorphism, i.e. that for a based quasimap
$\varphi\colon (C,\infty)\to(\CB,B_-)$ without defect at $0\in C$ with
$\varphi(0)\in U\cdot\{B_-\}\subset\CB$ there exists a unique
$\phi^C\in\dQM^\alpha_\fg$ such that $\varphi=(\phi^C(\infty))^{-1}\phi^C$.

Let $\varphi(0)=n'\cdot\{B_-\}$ for some $n'\in U$. We are looking for the
desired $\phi^C$ in the form $n^{-1}\varphi,\ n\in U$. So we must have
$\on{Ad}_{t_0}(n)=n^{-1}n'$, that is $t_0n^{-1}t_0^{-1}=n^{-1}n'\
\Leftrightarrow\ [n,t_0]=n'$. It remains to recall the following well known

\lem{comm} Let $t_0\in T^\reg$. Then the commutator map
$U\to U,\ n\mapsto[n,t_0]$ is an isomorphism (of algebraic varieties).
\elem

\prf
Filtering $U$ by its lower central series, one can introduce a system of
coordinates $(x_{i,j})_{1\leq i\leq h-1}^{1\leq j\leq b_i}$ on the affine space
$U$ such that for the inversion morphism $U\to U,\ n\mapsto n^{-1}$ we have
$(x_{i,j})^{-1}=(y_{i,j}(\unl{x}))$, and
$y_{i,j}(\unl{x})=-x_{i,j}+P_{i,j}(x_{i',j'})_{1\leq i'<i}^{1\leq j'\leq b_{i'}}$ for a certain polynomial $P_{i,j}$. Moreover, for the multiplication morphism
$m\colon U\times U\to U$ we have
$m((x'_{i,j}),(x''_{i,j}))=(x_{i,j}(\unl{x'},\unl{x''}))$, and
$x_{i,j}(\unl{x'},\unl{x''})=x'_{i,j}+x''_{i,j}+
Q_{i,j}(x'_{i',j'},x''_{i'',j''})_{1\leq i',i''<i}^{1\leq j'\leq b_{i'},
1\leq j''\leq b_{i''}}$ for a certain polynomial $Q_{i,j}$.
Finally, for the adjoint action $\on{Ad}_{t_0}\colon U\to U$ we have
$\on{Ad}_{t_0}(x_{i,j})=(w_{i,j}(\unl{x}))$, and
$w_{i,j}(\unl{x})=a_{i,j}x_{i,j}$ for a certain number $a_{i,j}\ne1$
(due to the regularity assumption on $t_0$).

Now given $n'=(x'_{i,j})\in U$ we can construct a unique $n=(x_{i,j}) \in U$ such that $[n,t_0]=n'$ recursively starting from $i=1$, and going to $i=h-1$.
\epr

The proposition is proved.
\epr

\cor{dense}
$Y^\alpha$ is an irreducible scheme with an open dense smooth subscheme
$\oY^\alpha$.
\ecor

\prf
We have seen that $q\colon Y^\alpha\to T^\reg$ is a fibration with a smooth
irreducible base and a fiber $F$ that is a $T$-torsor over $\hZ^\alpha$.
Now $\hZ^\alpha$ is open in the irreducible zastava scheme $Z^\alpha$ possessing
an open dense smooth subscheme $\oZ^\alpha$. Finally, 
$q\colon \oY^\alpha\to T^\reg$ is a fibration with a smooth irreducible base
and a fiber $\overset{\circ}{F}$ that is a $T$-torsor over 
$\dZ^\alpha=\oZ^\alpha\cap\hZ^\alpha$.
\epr

\sec{sympl}{A trigonometric symplectic structure}

\ssec{coord}{Coordinates on $Y^\alpha$}
Recall the locally closed embedding $Y^\alpha\hookrightarrow\hQM^\alpha_\fg$
introduced in the proof of~\refp{repr}. Via the Pl\"ucker embedding,
$\hQM^\alpha_\fg$ is a locally closed subscheme in
$\prod_{i\in I}V_{\check\omega_i}\otimes\Gamma(C,\CO_C(a_i))$ (notations
of~\refss{BunT}) where $\alpha=\sum_{i\in I}a_i\alpha_i$. In particular,
we have the coefficients $Q_i,R_i,S_{ij}\in\Gamma(C,\CO_C(a_i)) $ of the
highest, prehighest and next highest vectors
$v_{\check\omega_i},F_iv_{\check\omega_i},F_jF_iv_{\o_i}$.
Thus $Q_i,R_i,S_{ij}$ are the regular functions on $Y^\alpha$ with coefficients
in the space of degree $\leq a_i$ polynomials in $z$. The conditions
in~\refss{Yalpha}(c) ensure that $\deg Q_i=a_i$, and $Q_i(0)\ne0$.

It follows from~\cite[Remark 2]{fkmm} and~\refp{diamond} that (the coefficients
of) $(Q_i,R_i)_{i\in I}$ form a rational coordinate system on $Y^\alpha$.
Let us denote by $B_i$ (resp. $b_i$) the leading coefficient (resp. constant
term) of $Q_i$, so that $Q_i=B_iz^{a_i}+\ldots+b_i$. Similarly, we have
$R_i=C_iz^{a_i}+\ldots+c_i$. Note that $B_i\ne0\ne b_i$.
Following~\cite[3.3]{fkmm}, we introduce a rational \'etale
coordinate system on $Y^\alpha$. Namely,
$(w_{i,r})_{i\in I}^{1\leq r\leq a_i}$ are the ordered roots of $Q_i$, and
$y_{i,r}:=B_i^{-1}
R_i(w_{i,r})$. The desired coordinate system is formed by
$(B_i,C_i,w_{i,r},y_{i,r})_{i\in I}^{1\leq r\leq a_i}$.
It follows from~\cite[Remark 2]{fkmm} and~\refp{diamond} that these functions
do form a coordinate system on an unramified covering of the open subset of
$\oY^\alpha$ where all the roots of all the polynomials $Q_i,\ i\in I$, are
distinct.

We describe the $T\times T$-action on $Y^\alpha$, and the morphisms
$(p,q)\colon Y^\alpha\to T\times T$ (see~\refss{Yalpha}) in the above coordinates.
Note that the collection of fundamental weights $\check\omega_i\colon T\to\BG_m$
identifies $T$ with $\BG_m^I$. We have $(t_1,t_2)\cdot(Q_i,R_i)=
(\check\omega_i(t_1t_2)Q_i,\check\omega_i(t_1t_2)\check\alpha_i(t_1)^{-1}R_i)$,
and $\check\omega_j(p(Q_i,R_i)_{i\in I})=B_j^{-1}b_j$, and
$\check\alpha_j(q(Q_i,R_i)_{i\in I})=B_j^{-1}c_j^{-1}b_jC_j$.

\ssec{tang}{The tangent bundle}
Our goal in this Section is to describe the tangent space $T_y\oY^\alpha$ at
$y=(\CF_T,\tau_\fc,\phi,f_\fc)\in\oY^\alpha$. We denote by $\fg,\fb,\fu,\ft$
the Lie algebras of $G,B,U,T$. Given a $T$-bundle $\CF_T$
we denote the vector bundle associated to the adjoint action of $T$ on $\fg$ by
$\fg^\CF$. It is a direct sum of two subbundles corresponding to the trivial
(resp. nontrivial) eigenvalues of $T$ on $\fg\colon \fg^\CF=\ft^\CF\oplus\frr^\CF$.
Note that $\ft^\CF=\ft\otimes\CO_\Cx$. A $B$-structure $\phi$ on $\CF_G$ gives
rise to a vector subbundle $\fb^\phi\subset\fg^\CF$. The adjoint action of
$B$ on $\fu$ gives rise to a subbundle $\fu^\phi\subset\fb^\phi$. We denote
the quotient bundle by $\fh^\phi=\fh\otimes\CO_\Cx$: a trivial bundle where
$\fh=\fb/\fu$ is the abstract Cartan. The Killing form identifies the dual
vector bundle $\fb^{\phi*}$ with the quotient bundle
$\fg^\CF/\fu^\phi=:(\fg/\fu)^\phi$. For a vector bundle $\CV$ on $\Cx$ we
denote by $\CV_\fc$ the skyscraper quotient of $\CV$ by the ideal sheaf of the
point $\fc\in\Cx$.

We consider the following complex $K_y^\bullet$ of coherent sheaves on $\Cx$:
it lives in degrees $-1,0$, and $K_y^{-1}=(\ft\oplus\fh)\otimes\pi_*\omega_C$,
while $K_y^0=(\fg/\fu)^\phi$. The differential $d\colon K_y^{-1}\to K_y^0$ is a
direct sum of $d'\colon \ft\otimes\pi_*\omega_C\to(\fg/\fu)^\phi$ and
$d''\colon \fh\otimes\pi_*\omega_C\to(\fg/\fu)^\phi$. Here $d'$ is the
composition of $\ft\otimes\pi_*\omega_C\hookrightarrow\ft\otimes\CO_\Cx=
\ft^\CF$ (see~\refss{Cdag}) and
$\ft^\CF\hookrightarrow\fg^\CF\twoheadrightarrow(\fg/\fu)^\phi$, while
$d''$ is the composition  $\fh\otimes\pi_*\omega_C\hookrightarrow\fh\otimes\CO_\Cx=
\fh^\phi\hookrightarrow(\fg/\fu)^\phi$.

\prop{tange}
There is a canonical isomorphism $T_y\oY^\alpha\cong H^0(\Cx,K_y^\bullet)$.
\eprop

\prf
We consider the following complex $'\!K_y^\bullet$ of coherent sheaves on $\Cx$:
it lives in degrees $-1,0$, and $'\!K_y^{-1}=\fb^\phi$, while
$'\!K_y^0=\ft^\CF_\fc\oplus\fh^\phi_\fc$. The differential from $'\!K_y^{-1}$ to
$'\!K_y^0$ is a direct sum of $d'\colon \fb^\phi\to\ft^\CF_\fc$ and
$d''\colon \fb^\phi\to\fh^\phi_\fc$ where $d'$ is the composition
$\fb^\phi\hookrightarrow\fg^\CF\twoheadrightarrow\ft^\CF\twoheadrightarrow
\ft^\CF_\fc$, and $d''$ is the composition
$\fb^\phi\twoheadrightarrow\fh^\phi\twoheadrightarrow\fh^\phi_\fc$.

Then $T_y\oY^\alpha=H^0(\Cx,\ '\!K_y^\bullet)$. Now consider yet another complex
$''\!K_y^\bullet$ living in degrees $-1,0$ such that $''\!K_y^{-1}=\ '\!K_y^{-1}$
and $''\!K_y^0=\ '\!K_y^0\oplus\frr^\CF$, and the differential equals
$d'+d''+d'''$ where $d'''\colon \fb^\phi\to\frr^\CF$ is the composition
$\fb^\phi\hookrightarrow\fg^\CF\twoheadrightarrow\frr^\CF$. We have a canonical
morphism $''\!K_y^\bullet\to\ '\!K_y^\bullet$ inducing an isomorphism on
cohomology $H^0(\Cx,\ ''\!K_y^\bullet)\iso H^0(\Cx,\ '\!K_y^\bullet)$ since
$H^\bullet(\Cx,\frr^\CF)=0$ due to the regularity assumption on $\CF_T$.

Also we have a canonical quasiisomorphism $'''\!K_y^\bullet\to\ ''\!K_y^\bullet$
where $'''\!K_y^{-1}=\ ''\!K_y^{-1}\oplus(\ft\oplus\fh)\otimes\pi_*\omega_C$,
and $'''\!K_y^0=(\ft\oplus\fh)\otimes\CO_\Cx\oplus\frr^\CF$. The new components
of the differential are as follows:
$(\ft\oplus\fh)\otimes\pi_*\omega_C\hookrightarrow(\ft\oplus\fh)\otimes\CO_\Cx$
(see~\refss{Cdag}), and $\fb^\phi\twoheadrightarrow\fh^\phi$, and the composition $\fb^\phi\hookrightarrow\fg^\CF\twoheadrightarrow\ft^\CF$.

Finally, note that $'''\!K_y^0=\fh^\phi\oplus\fg^\CF$, and
$(\fh^\phi\oplus\fg^\CF)/d(\fb^\phi)\cong(\fg/\fu)^\phi$, so we have a
canonical quasiisomorphism $'''\!K_y^\bullet\to K_y^\bullet$.

The composition of the morphisms induced on $H^0(\Cx,?)$ by the above
quasiisomorphisms is the desired isomorphism
$T_y\oY^\alpha=H^0(\Cx,\ '\!K_y^\bullet)\cong H^0(\Cx,K_y^\bullet)$.
\epr

\ssec{cotang}{The cotangent bundle}
Let us describe the Serre dual complex $L_y^\bullet:=\CalD K_y^\bullet$.
It lives in degrees $0,1$, and $L_y^0=\fb^\phi$, while
$L_y^1=(\ft\oplus\fh)\otimes\pi_*\CO_C$.
The differential $d\colon L_y^0\to L_y^1$ is a direct sum of
$d'\colon \fb^\phi\to\ft\otimes\pi_*\CO_C$ and
$d''\colon \fb^\phi\to\fh\otimes\pi_*\CO_C$. Here $d''$ is the composition
$\fb^\phi\twoheadrightarrow\fh^\phi=\fh\otimes\CO_\Cx\hookrightarrow
\fh\otimes\pi_*\CO_C$ (see~\refss{Cdag}), while $d'$ is the composition
$\fb^\phi\hookrightarrow\fg^\CF\twoheadrightarrow\ft^\CF=\ft\otimes\CO_\Cx
\hookrightarrow\ft\otimes\pi_*\CO_C$. Now~\refp{tange} has the following
immediate

\cor{cotange}
There is a canonical isomorphism $T^*_y\oY^\alpha\cong H^1(\Cx,L_y^\bullet)$.
\qed
\ecor

\ssec{somedi}{Some differentials}
Here we describe the differentials of the morphism
$(p,q)\colon \oY^\alpha\to T\times T^\reg$ and of the action of $T\times T$
on $\oY^\alpha$ introduced in~\refss{Yalpha}.

Note that for a regular $T$-bundle $\CF_T$ on $\Cx$ the tangent space
$T_{\CF_T}\Bun_T^\reg(\Cx)$ is canonically isomorphic to
$H^1(\Cx,\ft^\CF)=H^1(\Cx,\ft\otimes\CO_\Cx)\cong\ft$ (here the second
isomorhism is $\on{Id}_\ft\otimes\on{Tr}$ for the trace isomorphism
$\on{Tr}\colon H^1(\Cx,\CO_\Cx)=H^1(\Cx,\omega_\Cx)\iso\BC$).
The scalar product $(,)$
on $\ft^*$~(\refss{BunT}) identifies $\ft$ with $\ft^*$. Together with the
Serre duality $H^1(\Cx,\CO_\Cx)^*=H^0(\Cx,\CO_\Cx)$ this gives rise to a
canonical isomorphism
$T^*_{\CF_T}\Bun_T^\reg(\Cx)\cong H^0(\Cx,\ft\otimes\CO_\Cx)\cong\ft$.
Similarly, for a degree $\alpha$ $T$-bundle $\phi_T$ we have canonical
isomorphisms
$T_{\phi_T}\Bun^\alpha_T(\Cx)\cong H^1(\Cx,\fh^\phi)=H^1(\fh\otimes\CO_\Cx)=\fh$
and $T^*_{\phi_T}\Bun_T^\alpha(\Cx)\cong H^0(\Cx,\fh\otimes\CO_\Cx)=\fh$.

We have distinguished triangles
$(\fg/\fu)^\phi\to K^\bullet_y\to(\ft\oplus\fh)\otimes\pi_*\omega_C[1]$ and
$(\ft\oplus\fh)\otimes\pi_*\CO_C[-1]\to L^\bullet_y\to\fb^\phi$.
They give rise to a morphism $$\fd_p\colon T_y\oY^\alpha=H^0(\Cx,K^\bullet_y)\to
H^1(\Cx,(\ft\oplus\fh)\otimes\pi_*\omega_C)\to$$
$$\to H^1(\Cx,(\ft\oplus\fh)\otimes\omega_\Cx)\to H^1(\Cx,\fh\otimes\omega_\Cx)
=T_{\phi_T}\Bun^\alpha_T(\Cx)$$ where the middle arrow arises from the natural
morphism $\pi_*\omega_C\to\omega_\Cx$ (see~\refss{Cdag}), and the next
arrow arises from the projection $\ft\oplus\fh\to\fh$. Similarly, we have
$$\fd_q\colon T_y\oY^\alpha=H^0(\Cx,K^\bullet_y)\to
H^1(\Cx,(\ft\oplus\fh)\otimes\pi_*\omega_C)\to$$
$$\to H^1(\Cx,(\ft\oplus\fh)\otimes\omega_\Cx)\to H^1(\Cx,\ft\otimes\omega_\Cx)
=T_{\CF_T}\Bun^\reg_T(\Cx).$$ Dually, we have 
$$\fd_p^*\colon T^*_{\phi_T}\Bun^\alpha_T(\Cx)
=H^0(\Cx,\fh\otimes\CO_\Cx)\to H^0(\Cx,(\ft\oplus\fh)\otimes\CO_\Cx)\to$$
$$\to H^0(\Cx,(\ft\oplus\fh)\otimes\pi_*\CO_C)\to H^1(\Cx,L^\bullet_y)=
T_y^*\oY^\alpha$$ and $$\fd_q^*\colon T^*_{\CF_T}\Bun^\reg_T(\Cx)
=H^0(\Cx,\ft\otimes\CO_\Cx)\to H^0(\Cx,(\ft\oplus\fh)\otimes\CO_\Cx)\to$$
$$\to H^0(\Cx,(\ft\oplus\fh)\otimes\pi_*\CO_C)\to H^1(\Cx,L^\bullet_y)=
T_y^*\oY^\alpha.$$

We also have a morphism $$\fa_1\colon T_eT=\ft=H^0(\Cx,\ft\otimes\CO_\Cx)=
H^0(\Cx,\ft^\CF)\to H^0(\Cx,(\fg/\fu)^\phi)\to H^0(\Cx,K^\bullet_y)=T_y\oY^\alpha,$$ 
and
$$\fa_2\colon T_eT\cong\fh=H^0(\Cx,\fh\otimes\CO_\Cx)=
H^0(\Cx,\fh^\phi)\to H^0(\Cx,(\fg/\fu)^\phi)\to H^0(\Cx,K^\bullet_y)=T_y\oY^\alpha.$$

\lem{somedif}
a) $(\fa_1,\fa_2)\colon \ft\oplus\fh\to T_y\oY^\alpha$ is the differential of the
action $T\times T\times\oY^\alpha\to\oY^\alpha$ (see~\refss{Yalpha});

b) $(\fd_p,\fd_q)\colon T_y\oY^\alpha\to T_{\phi_T}\Bun^\alpha_T(\Cx)\oplus
T_{\CF_T}\Bun^\reg_T(\Cx)$ is the differential of
$(p,q)\colon \oY^\alpha\to\Bun^\alpha_T(\Cx)\times\Bun^\reg_T(\Cx)$.

c) $(\fd_p^*,\fd_q^*)\colon T_{\phi_T}^*\Bun^\alpha_T(\Cx)\oplus
T_{\CF_T}^*\Bun^\reg_T(\Cx)\to T_y^*\oY^\alpha$ is the codifferential of
$(p,q)\colon \oY^\alpha\to\Bun^\alpha_T(\Cx)\times\Bun^\reg_T(\Cx)$.
\elem

\prf Clear from the construction.
\epr

\ssec{bivec}{A bivector field}
We consider the following bicomplex $M_y^{\bullet,\bullet}$:
\eq{bico}
\begin{CD}
 @. \fb^\phi @>>> (\ft\oplus\fh)\otimes\pi_*\CO_C \\
@.    @VVV                                            @VVV \\
(\ft\oplus\fh)\otimes\pi_*\omega_C @>>> \frr^\CF\oplus(\ft\oplus\fh)\otimes
(\theta\oplus\CO_\Cx) @>>> (\ft\oplus\fh)\otimes\pi_*\CO_C  \\
@VVV                                 @VVV    @. \\
(\ft\oplus\fh)\otimes\pi_*\omega_C @>>> (\fg/\fu)^\phi @. \\
\end{CD}
\end{equation}
Here the middle term lives in bidegree $(0,0)$, the first line is nothing but
$L_y^\bullet$ of~\refss{cotang}, while the last line is nothing but
$K_y^\bullet$ of~\refss{tang}. The left vertical arrow is the identity morphism,
as well as the right vertical arrow. The middle line is a direct sum of the
complex consisting of $\frr^\CF$ in degree 0, and of the exact
complex~\refe{theta omega} tensored with $\ft\oplus\fh$. Note that the middle
term can be rearranged as
$\fh^\phi\oplus\fg^\CF\oplus(\ft\oplus\fh)\otimes\theta$. Now the middle column
is a direct sum of the complex consisting of $(\ft\oplus\fh)\otimes\theta$
in degree 0, and of the exact complex
$\fb^\phi\to\fh^\phi\oplus\fg^\CF\to(\fg/\fu)^\phi$.

It follows (looking at the columns of~\refe{bico}) that the total complex
$\on{Tot}M_y^{\bullet,\bullet}$ has only one cohomology in degree 0, and
$H^0(\on{Tot}M_y^{\bullet,\bullet})=(\ft\oplus\fh)\otimes\theta$.
Since $H^\bullet(\Cx,(\ft\oplus\fh)\otimes\theta)=0$, we deduce that
$H^\bullet(\Cx,\on{Tot}M_y^{\bullet,\bullet})=0$. Now let us look at the rows
of $M_y^{\bullet,\bullet}$. The hypercohomology of $\Cx$ with coefficients in
the middle row vanishes since $H^\bullet(\Cx,\frr^\CF)=0$ (due to the regularity
assumption on $\CF_T$). Hence the second differential in the spectral sequence
converging to $H^\bullet(\Cx,\on{Tot}M_y^{\bullet,\bullet})=0$ from the
hypercohomology of $\Cx$ with coefficients in the rows is
$d_2\colon H^1(\Cx,L_y^\bullet)\iso H^0(\Cx,K_y^\bullet)$.

Finally, due to~\refp{tange} and~\refc{cotange} we can view the above
differential as $d_2\colon T_y^*\oY^\alpha\iso T_y\oY^\alpha$.

%\lem{mom}
%The following diagrams commute:\\
%$\begin{CD}
%\ft @= T_eT\cong \fh @>{\fa_2}>> T_y\oY^\alpha \\
%@|  @. @A{d_2}AA \\
%\ft^* @= T^*_{\CF_T}\Bun^\reg_T(\Cx) @>{\fd^*_q}>> T^*_y\oY^\alpha\\
%\end{CD}$ and
%$\begin{CD}
%\fh @= \ft\cong T_eT @>{\fa_1}>> T_y\oY^\alpha \\
%@|  @. @A{d_2}AA \\
%\fh^* @= T^*_{\phi_T}\Bun^\alpha_T(\Cx) @>{\fd^*_p}>> T^*_y\oY^\alpha\\
%\end{CD}$\\
%(the left vertical equalities are the identifications via the form $(,)$
%of~\refss{BunT}).
%\elem

%\prf
%Easy but wrong.
%\epr

\ssec{calcul}{Calculation of the bivector field: preparation}
We follow the strategy of~\cite{fkmm}, and eventually reduce our calculation to
the one of {\em loc. cit.} The result of the somewhat lengthy calculation of
$d_2$ is contained in~\refl{quad},~\refr{discla} and~\refp{atlast}.

Given a character $\check\lambda\in\Lambda^\vee$, we consider the composed
homomorphism $B\to T\to\BG_m$, and denote the associated (to the $B$-torsor
$\phi$) line bundle on $\Cx$ by $\CL_{\check\lambda}^\phi$.
For an irreducible $G$-module $V_{\check\lambda}$, the associated
(to the $G$-torsor $\CF_G$) vector bundle
on $\Cx$ is denoted $\CV^\CF_{\check\lambda}$. If $\check\lambda$ is a
fundamental weight $\o_i$, then we have an isomorphism
$V^*_{\o_i}\cong V_{\o_{i^*}}$ for an involution $I\iso I,\ i\mapsto i^*$.
If we extend the involution $\o_i\mapsto\o_{i^*}$ by linearity to the weight
lattice $\Lambda^\vee,\ \check\lambda\mapsto\check\lambda^*$, then this
involution preserves the scalar product $(,)$ of~\refss{BunT}.

We have a natural embedding of vector bundles on
$\Cx\colon (\fg/\fu)^\phi\hookrightarrow\bigoplus_{i\in I}\CV^\CF_{\o_{i^*}}\otimes
\CL_{\o_i}^\phi$, and the dual surjection $\bigoplus_{i\in I}
\CV^{*\CF}_{\o_{i^*}}\otimes\CL_{-\o_i}^\phi\twoheadrightarrow\fb^\phi$.
They give rise to the following morphisms of two-term complexes of coherent
sheaves on $\Cx$:
\eq{bottom}
\begin{CD}
(\ft\oplus\fh)\otimes\pi_*\omega_C @>>> (\fg/\fu)^\phi \\
@VVV                                       @VVV \\
(\ft\oplus\fh)\otimes\pi_*\Xi%\theta
@>>> (\fg/\fu)^\phi\otimes\pi_*\CO_C \\
@|                                       @VVV\\
(\ft\oplus\fh)\otimes\pi_*\Xi%\theta
@>>> (\bigoplus_{i\in I}\CV^\CF_{\o_{i^*}}\otimes
\CL_{\o_i}^\phi)\otimes\pi_*\CO_C\\
\end{CD}
\end{equation}
(the upper vertical arrows arise from the morphisms $\pi_*\omega_C\to\pi_*\Xi=
\pi_*(\CO_C(-\{0\})\oplus\CO_C(-\{\infty\}))$
and $\CO_\Cx\to\pi_*\CO_C$ of~\refss{Cdag}), and dually
\eq{top}
\begin{CD}
(\bigoplus_{i\in I}
\CV^{*\CF}_{\o_{i^*}}\otimes\CL_{-\o_i}^\phi)\otimes\pi_*\omega_C @>>>
(\ft\oplus\fh)\otimes\pi_*\Xi%\theta
\\
@VVV       @| \\
\fb^\phi\otimes\pi_*\omega_C @>>>
(\ft\oplus\fh)\otimes\pi_*\Xi%\theta
\\
@VVV                                   @VVV \\
\fb^\phi @>>> (\ft\oplus\fh)\otimes\pi_*\CO_C \\
\end{CD}
\end{equation}
Note that the top row of~\refe{bottom} coincides with the bottom row
$K^\bullet_y$ of the bicomplex~\refe{bico}, while the bottom row
of~\refe{top} coincides with the top row $L^\bullet_y$ of the
bicomplex~\refe{bico}. So composing the vertical arrows
of~\refe{top},~\refe{bottom} with the vertical arrows of~\refe{bico} we
obtain a bicomplex $'\!M_y^{\bullet,\bullet}$:
\eq{bicompl}
\begin{CD}
 @. (\bigoplus_{i\in I}
\CV^{*\CF}_{\o_{i^*}}\otimes\CL_{-\o_i}^\phi)\otimes\pi_*\omega_C @>>>
(\ft\oplus\fh)\otimes\pi_*\Xi%\theta
\\
@.    @VVV                                            @VVV \\
(\ft\oplus\fh)\otimes\pi_*\omega_C @>>> \frr^\CF\oplus(\ft\oplus\fh)\otimes
(\theta\oplus\CO_\Cx) @>>> (\ft\oplus\fh)\otimes\pi_*\CO_C  \\
@VVV                                 @VVV    @. \\
(\ft\oplus\fh)\otimes\pi_*\Xi%\theta
@>>> (\bigoplus_{i\in I}\CV^\CF_{\o_{i^*}}\otimes
\CL_{\o_i}^\phi)\otimes\pi_*\CO_C @. \\
\end{CD}
\end{equation}
Note that
$H^\bullet(\Cx,(\ft\oplus\fh)\otimes\pi_*\Xi)=0$, so just as
in~\refss{bivec}, the second differential in the spectral sequence converging
to $H^\bullet(\Cx,\on{Tot}\,'\!M_y^{\bullet,\bullet})$ from the hypercohomology of $\Cx$
with coefficients in the rows is
\eq{d'}
d'_2\colon H^1(\Cx,(\bigoplus_{i\in I}
\CV^{*\CF}_{\o_{i^*}}\otimes\CL_{-\o_i}^\phi)\otimes\pi_*\omega_C)\to H^0(\Cx,
(\bigoplus_{i\in I}\CV^\CF_{\o_{i^*}}\otimes
\CL_{\o_i}^\phi)\otimes\pi_*\CO_C).
\end{equation}

\lem{quad}
The following diagram commutes:

$\begin{CD}
H^0(\Cx,L^\bullet_y) @<<< H^1(\Cx,(\bigoplus_{i\in I}
\CV^{*\CF}_{\o_{i^*}}\otimes\CL_{-\o_i}^\phi)\otimes\pi_*\omega_C) \\
@V{d_2}VV                @V{d'_2}VV \\
H^0(\Cx,K^\bullet_y) @>>> H^0(\Cx,(\bigoplus_{i\in I}\CV^\CF_{\o_{i^*}}\otimes
\CL_{\o_i}^\phi)\otimes\pi_*\CO_C) \\
\end{CD}$
\elem

\prf
Clear.
\epr

\rem{discla}
In what follows we will be occupied with the calculation of
$$d'_2\colon H^1(\Cx,(\bigoplus_{i\in I}
\CV^{*\CF}_{\o_{i^*}}\otimes\CL_{-\o_i}^\phi)\otimes\pi_*\omega_C)\to H^0(\Cx,
(\bigoplus_{i\in I}\CV^\CF_{\o_{i^*}}\otimes
\CL_{\o_i}^\phi)\otimes\pi_*\CO_C).$$ Let us presently comment in which sense
does it calculate the desired
$d_2\colon H^1(\Cx,L^\bullet_y)\to H^0(\Cx,K^\bullet_y)$.
It is easy to see that for $y\in\oY^\alpha$ lying in the open subset
$U^\alpha\subset\oY^\alpha$ formed by all the quadruples $(\CF_T,\tau_\fc,\phi,f_\fc)$
such that $\ft^\CF_\fc\cap\fb^\phi_\fc=0\subset\fg^\CF_\fc$ the morphism
$H^1(\Cx,(\bigoplus_{i\in I}
\CV^{*\CF}_{\o_{i^*}}\otimes\CL_{-\o_i}^\phi)\otimes\pi_*\omega_C)\to
H^1(\Cx,L^\bullet_y)$ is surjective, and the morphism
$H^0(\Cx,K^\bullet_y)\to H^0(\Cx,(\bigoplus_{i\in I}\CV^\CF_{\o_{i^*}}\otimes
\CL_{\o_i}^\phi)\otimes\pi_*\CO_C)$ is injective. Since we are only going to
calculate our $d_2$ generically, the only trouble is that for some $\alpha$
the open subset $U^\alpha$ may happen to be empty. Indeed, for
$\alpha=\sum_{i\in I}a_i\alpha_i$, we have $U^\alpha=\emptyset$ iff
$a_i=0$ for some $i\in I$. So in what follows we assume $a_i>0\ \forall i\in I$
(otherwise the moduli space $\oY^\alpha$ essentially reduces to the one for
a semisimple Lie algebra $\fg'$ of smaller rank).
\erem

\ssec{reducto}{Reduction to a calculation on $C$}
%Recall from~\refss{bivec} and~\refss{calcul} that the middle column of the
%bicomplex $'\!M_y^{\bullet,\bullet}$~\refe{bicompl} is the direct sum of
%the complex consisting of $(\ft\oplus\fh)\otimes\theta$ in degree 0, and of
%the complex
%\eq{column}
%(\bigoplus_{i\in I}
%\CV^{*\CF}_{\o_{i^*}}\otimes\CL_{-\o_i}^\phi)\otimes\pi_*\omega_C\to
%\fh^\phi\oplus\fg^\CF\to(\bigoplus_{i\in I}\CV^\CF_{\o_{i^*}}\otimes
%\CL_{\o_i}^\phi)\otimes\pi_*\CO_C.
%\end{equation}
%The second differential in the spectral
%sequence converging to the hypercohomology of $\Cx$ with coefficients in the
%latter complex is $d_2\colon H^1(\Cx,(\bigoplus_{i\in I}
%\CV^{*\CF}_{\o_{i^*}}\otimes\CL_{-\o_i}^\phi)\otimes\pi_*\omega_C)\to H^0(\Cx,
%(\bigoplus_{i\in I}\CV^\CF_{\o_{i^*}}\otimes
%\CL_{\o_i}^\phi)\otimes\pi_*\CO_C)$. Clearly, $d'_2=d_2$.
The goal of this Section is a description of $d'_2$~\refe{d'} in terms of
$C$, see~\refc{botline}.

Note that $\pi^*(\fh^\phi\oplus\fg^\CF)\cong(\fh\oplus\fg)\otimes\CO_C$,
and $\pi^*((\ft\oplus\fh)\otimes\theta)\cong(\ft\oplus\fh)\otimes\CO_C$.
Hence $\pi_*\pi^*(\fh^\phi\oplus\fg^\CF\oplus(\ft\oplus\fh)\otimes\theta)\cong
(\fh\oplus\fg\oplus\ft\oplus\fh)\otimes\pi_*\CO_C$.
The morphisms $$\frr^\CF\oplus(\ft\oplus\fh)\otimes
(\theta\oplus\CO_\Cx)=\fh^\phi\oplus\fg^\CF\oplus(\ft\oplus\fh)\otimes\theta\to (\ft\oplus\fh)\otimes\pi_*\CO_C,$$ 
$$\on{resp.}\ \frr^\CF\oplus(\ft\oplus\fh)\otimes
(\theta\oplus\CO_\Cx)=\fh^\phi\oplus\fg^\CF\oplus(\ft\oplus\fh)\otimes\theta\to
(\bigoplus_{i\in I}\CV^\CF_{\o_{i^*}}\otimes
\CL_{\o_i}^\phi)\otimes\pi_*\CO_C$$ of~\refe{bicompl} factor through the
canonical morphism $$\fh^\phi\oplus\fg^\CF\oplus(\ft\oplus\fh)\otimes\theta\to
\pi_*\pi^*(\fh^\phi\oplus\fg^\CF\oplus(\ft\oplus\fh)\otimes\theta)\cong
(\fh\oplus\fg\oplus\ft\oplus\fh)\otimes\pi_*\CO_C$$ and
$$(\fh\oplus\fg\oplus\ft\oplus\fh)\otimes\pi_*\CO_C\to
(\ft\oplus\fh)\otimes\pi_*\CO_C,$$ 
$$\on{resp.}\ (\fh\oplus\fg\oplus\ft\oplus\fh)\otimes\pi_*\CO_C\to
(\bigoplus_{i\in I}\CV^\CF_{\o_{i^*}}\otimes
\CL_{\o_i}^\phi)\otimes\pi_*\CO_C.$$ Hence we obtain a morphism from the
bicomplex $'\!M_y^{\bullet,\bullet}$~\refe{bicompl} to the following bicomplex
$''\!M_y^{\bullet,\bullet}$:
\eq{bicomple}
\begin{CD}
 @. (\bigoplus_{i\in I}
\CV^{*\CF}_{\o_{i^*}}\otimes\CL_{-\o_i}^\phi)\otimes\pi_*\omega_C @>>>
(\ft\oplus\fh)\otimes\pi_*\Xi%\theta
\\
@.    @VVV                                            @VVV \\
(\ft\oplus\fh)\otimes\pi_*\omega_C @>>>
(\fh\oplus\fg\oplus\ft\oplus\fh)\otimes\pi_*\CO_C
@>>> (\ft\oplus\fh)\otimes\pi_*\CO_C  \\
@VVV                                 @VVV    @. \\
(\ft\oplus\fh)\otimes\pi_*\Xi%\theta
@>>> (\bigoplus_{i\in I}\CV^\CF_{\o_{i^*}}\otimes
\CL_{\o_i}^\phi)\otimes\pi_*\CO_C @. \\
\end{CD}
\end{equation}
(the morphisms from all the terms of $'\!M_y^{\bullet,\bullet}$ to the
corresponding terms of $''\!M_y^{\bullet,\bullet}$ except for the middle ones are identities). Just as in~\refss{calcul} we obtain the second differential
\eq{d''}
d''_2\colon H^1(\Cx,(\bigoplus_{i\in I}
\CV^{*\CF}_{\o_{i^*}}\otimes\CL_{-\o_i}^\phi)\otimes\pi_*\omega_C)\to H^0(\Cx,
(\bigoplus_{i\in I}\CV^\CF_{\o_{i^*}}\otimes
\CL_{\o_i}^\phi)\otimes\pi_*\CO_C)
\end{equation}
in the spectral sequence converging
to $H^\bullet(\Cx,\on{Tot}\,''\!M_y^{\bullet,\bullet})$ from the hypercohomology of $\Cx$ with coefficients in the rows. It follows that $d''_2=d'_2$
of~\refe{d'}.

Now the morphisms $$(\bigoplus_{i\in I}
\CV^{*\CF}_{\o_{i^*}}\otimes\CL_{-\o_i}^\phi)\otimes\pi_*\omega_C\to
(\fh\oplus\fg\oplus\ft\oplus\fh)\otimes\pi_*\CO_C,$$ 
$$\on{resp.}\ (\ft\oplus\fh)\otimes\pi_*\omega_C\to
(\fh\oplus\fg\oplus\ft\oplus\fh)\otimes\pi_*\CO_C$$ of~\refe{bicomple}
factor through the natural morphism
$$(\fh\oplus\fg\oplus\ft\oplus\fh)\otimes\pi_*\Xi\to
(\fh\oplus\fg\oplus\ft\oplus\fh)\otimes\pi_*\CO_C$$ (see~\refss{Cdag}) and
$$(\bigoplus_{i\in I}
\CV^{*\CF}_{\o_{i^*}}\otimes\CL_{-\o_i}^\phi)\otimes\pi_*\omega_C\to
(\fh\oplus\fg\oplus\ft\oplus\fh)\otimes\pi_*\Xi,$$ 
$$\on{resp.}\ (\ft\oplus\fh)\otimes\pi_*\omega_C\to
(\fh\oplus\fg\oplus\ft\oplus\fh)\otimes\pi_*\Xi.$$
Hence we obtain a morphism to the
bicomplex $''\!M_y^{\bullet,\bullet}$~\refe{bicompl} to the following bicomplex
$'''\!M_y^{\bullet,\bullet}$:
\eq{bicomplex}
\begin{CD}
 @. (\bigoplus_{i\in I}
\CV^{*\CF}_{\o_{i^*}}\otimes\CL_{-\o_i}^\phi)\otimes\pi_*\omega_C @>>>
(\ft\oplus\fh)\otimes\pi_*\Xi%\theta
\\
@.    @VVV                                            @VVV \\
(\ft\oplus\fh)\otimes\pi_*\omega_C @>>>
(\fh\oplus\fg\oplus\ft\oplus\fh)\otimes\pi_*\Xi
@>>> (\ft\oplus\fh)\otimes\pi_*\CO_C  \\
@VVV                                 @VVV    @. \\
(\ft\oplus\fh)\otimes\pi_*\Xi%\theta
@>>> (\bigoplus_{i\in I}\CV^\CF_{\o_{i^*}}\otimes
\CL_{\o_i}^\phi)\otimes\pi_*\CO_C @. \\
\end{CD}
\end{equation}
(the morphisms from all the terms of $'''\!M_y^{\bullet,\bullet}$ to the
corresponding terms of $''\!M_y^{\bullet,\bullet}$ except for the middle ones are identities). Just as in~\refss{calcul} we obtain the second differential
\eq{d'''}
d'''_2\colon H^1(\Cx,(\bigoplus_{i\in I}
\CV^{*\CF}_{\o_{i^*}}\otimes\CL_{-\o_i}^\phi)\otimes\pi_*\omega_C)\to H^0(\Cx,
(\bigoplus_{i\in I}\CV^\CF_{\o_{i^*}}\otimes
\CL_{\o_i}^\phi)\otimes\pi_*\CO_C)
\end{equation}
in the spectral sequence converging
to $H^\bullet(\Cx,\on{Tot}\,'''\!M_y^{\bullet,\bullet})$ from the hypercohomology of $\Cx$ with coefficients in the rows. It follows that $d'''_2=d''_2$
of~\refe{d''}.

Note that the bicomplex $'''\!M_y^{\bullet,\bullet}$ is obtained from
the following bicomplex $^\circ\!M_y^{\bullet,\bullet}$:
\eq{Nbicomplex}
\begin{CD}
 @. (\bigoplus_{i\in I}
\CV^{*\CF}_{\o_{i^*}}\otimes\CL_{-\o_i}^\phi)\otimes\pi_*\Xi @>>>
(\ft\oplus\fh)\otimes\pi_*\Xi%\theta
\\
@.    @VVV                                            @VVV \\
(\ft\oplus\fh)\otimes\pi_*\omega_C @>>>
(\fh\oplus\fg\oplus\ft\oplus\fh)\otimes\pi_*\Xi
@>>> (\ft\oplus\fh)\otimes\pi_*\CO_C  \\
@VVV                                 @VVV    @. \\
(\ft\oplus\fh)\otimes\pi_*\Xi%\theta
@>>> (\bigoplus_{i\in I}\CV^\CF_{\o_{i^*}}\otimes
\CL_{\o_i}^\phi)\otimes\pi_*\Xi @. \\
\end{CD}
\end{equation}
by composing the vertical arrows of~\refe{Nbicomplex} with the vertical
arrows of the following commutative diagrams:
\eq{tuda}
\begin{CD}
(\bigoplus_{i\in I}
\CV^{*\CF}_{\o_{i^*}}\otimes\CL_{-\o_i}^\phi)\otimes\pi_*\omega_C @>>>
(\ft\oplus\fh)\otimes\pi_*\Xi\\
@VVV @| \\
(\bigoplus_{i\in I}
\CV^{*\CF}_{\o_{i^*}}\otimes\CL_{-\o_i}^\phi)\otimes\pi_*\Xi @>>>
(\ft\oplus\fh)\otimes\pi_*\Xi\\
\end{CD}
\end{equation}
\eq{suda}
\begin{CD}
(\ft\oplus\fh)\otimes\pi_*\Xi%\theta
@>>> (\bigoplus_{i\in I}\CV^\CF_{\o_{i^*}}\otimes
\CL_{\o_i}^\phi)\otimes\pi_*\Xi\\
@| @VVV \\
(\ft\oplus\fh)\otimes\pi_*\Xi%\theta
@>>> (\bigoplus_{i\in I}\CV^\CF_{\o_{i^*}}\otimes
\CL_{\o_i}^\phi)\otimes\pi_*\CO_C\\
\end{CD}
\end{equation}
Similarly to~\refl{quad} we deduce the following
\lem{quadr}
The following diagram commutes:

$\begin{CD}
H^1(\Cx,(\bigoplus_{i\in I}
\CV^{*\CF}_{\o_{i^*}}\otimes\CL_{-\o_i}^\phi)\otimes\pi_*\Xi) @<<< H^1(\Cx,(\bigoplus_{i\in I}
\CV^{*\CF}_{\o_{i^*}}\otimes\CL_{-\o_i}^\phi)\otimes\pi_*\omega_C) \\
@V{d^\circ_2}VV                @V{d'''_2}VV \\
H^0(\Cx,(\bigoplus_{i\in I}\CV^\CF_{\o_{i^*}}\otimes
\CL_{\o_i}^\phi)\otimes\pi_*\Xi) @>>> H^0(\Cx,(\bigoplus_{i\in I}\CV^\CF_{\o_{i^*}}\otimes
\CL_{\o_i}^\phi)\otimes\pi_*\CO_C) \\
\end{CD}$

where $d^\circ_2\colon H^1(\Cx,(\bigoplus_{i\in I}
\CV^{*\CF}_{\o_{i^*}}\otimes\CL_{-\o_i}^\phi)\otimes\pi_*\Xi)\to H^0(\Cx,
(\bigoplus_{i\in I}\CV^\CF_{\o_{i^*}}\otimes
\CL_{\o_i}^\phi)\otimes\pi_*\Xi)$ is the second differential
in the spectral sequence converging
to $H^\bullet(\Cx,\on{Tot}\,^\circ\!M_y^{\bullet,\bullet})$ from the hypercohomology of $\Cx$ with coefficients in the rows. \qed
\elem

By the projection formula, $$(\bigoplus_{i\in I}\CV^\CF_{\o_{i^*}}\otimes
\CL_{\o_i}^\phi)\otimes\pi_*\CO_C\cong
\pi_*\pi^*(\bigoplus_{i\in I}\CV^\CF_{\o_{i^*}}\otimes
\CL_{\o_i}^\phi)=\pi_*\bigoplus_{i\in I}V_{\o_{i^*}}\otimes
\pi^*\CL_{\o_i}^\phi,$$ and similarly
$$(\bigoplus_{i\in I}\CV^\CF_{\o_{i^*}}\otimes
\CL_{\o_i}^\phi)\otimes\pi_*\Xi\cong\pi_*((\bigoplus_{i\in I}V_{\o_{i^*}}\otimes
\pi^*\CL_{\o_i}^\phi)\otimes\Xi),$$ 
$$(\bigoplus_{i\in I}
\CV^{*\CF}_{\o_{i^*}}\otimes\CL_{-\o_i}^\phi)\otimes\pi_*\omega_C\cong
\pi_*((\bigoplus_{i\in I}
V^*_{\o_{i^*}}\otimes\pi^*\CL_{-\o_i}^\phi)\otimes\omega_C),$$ 
$$(\bigoplus_{i\in I}
\CV^{*\CF}_{\o_{i^*}}\otimes\CL_{-\o_i}^\phi)\otimes\pi_*\Xi\cong
\pi_*((\bigoplus_{i\in I}V^*_{\o_{i^*}}
\otimes\pi^*\CL_{-\o_i}^\phi)\otimes\Xi).$$

Now since the upper-right and lower-left terms $(\ft\oplus\fh)\otimes\pi_*\Xi$
of the bicomplex~\refe{Nbicomplex} are acyclic sheaves on $\Cx$, the
differential $d^\circ_2$ coincides with the second differential from the
spectral sequence arising from the following {\em complex} on $\Cx$:
\eq{verthor}
(\bigoplus_{i\in I}
\CV^{*\CF}_{\o_{i^*}}\otimes\CL_{-\o_i}^\phi)\otimes\pi_*\Xi\to
(\fh\oplus\fg\oplus\ft\oplus\fh)\otimes\pi_*\Xi\to
(\bigoplus_{i\in I}\CV^\CF_{\o_{i^*}}\otimes
\CL_{\o_i}^\phi)\otimes\pi_*\Xi
\end{equation}
which is by construction a direct sum of the one term complex
$(\ft\oplus\fh)\otimes\pi_*\Xi$ in degree zero and the following complex:
\eq{verthor'}
(\bigoplus_{i\in I}
\CV^{*\CF}_{\o_{i^*}}\otimes\CL_{-\o_i}^\phi)\otimes\pi_*\Xi\to
(\fh\oplus\fg)\otimes\pi_*\Xi\to
(\bigoplus_{i\in I}\CV^\CF_{\o_{i^*}}\otimes
\CL_{\o_i}^\phi)\otimes\pi_*\Xi
\end{equation}
which in turn is nothing but the direct image $\pi_*N^\bullet_y$ of the
following complex $N^\bullet_y$ of vector bundles on $C$:
\eq{on C}
(\bigoplus_{i\in I}V^*_{\o_{i^*}}\otimes\pi^*\CL_{-\o_i}^\phi)\otimes\Xi\to
(\fh\oplus\fg)\otimes\Xi\to
(\bigoplus_{i\in I}V_{\o_{i^*}}\otimes
\pi^*\CL_{\o_i}^\phi)\otimes\Xi
\end{equation}

\lem{clear}
The following diagram commutes:

$\begin{CD}
H^1(\Cx,(\bigoplus_{i\in I}
\CV^{*\CF}_{\o_{i^*}}\otimes\CL_{-\o_i}^\phi)\otimes\pi_*\Xi) @=
H^1(C,(\bigoplus_{i\in I}V^*_{\o_{i^*}}\otimes
\pi^*\CL_{-\o_i}^\phi)\otimes\Xi) \\
@V{d^\circ_2}VV                @V{d^C_2}VV \\
H^0(\Cx,(\bigoplus_{i\in I}\CV^\CF_{\o_{i^*}}\otimes
\CL_{\o_i}^\phi)\otimes\pi_*\Xi) @=
H^0(C,(\bigoplus_{i\in I}V_{\o_{i^*}}\otimes
\pi^*\CL_{\o_i}^\phi)\otimes\Xi) \\
\end{CD}$

where $d^C_2\colon H^1(C,(\bigoplus_{i\in I}V_{\o_{i^*}}\otimes
\pi^*\CL_{\o_i}^\phi)\otimes\Xi)\to
H^0(C,(\bigoplus_{i\in I}V_{\o_{i^*}}\otimes
\pi^*\CL_{\o_i}^\phi)\otimes\Xi)$ is the second
differential in the spectral sequence converging
to $H^\bullet(C,N_y^\bullet)$ from the cohomology of $C$ with coefficients in
its terms. \qed
\elem

Recall that our goal is to compute the differential~\refe{d'}:
$$H^1(C,(\bigoplus_{i\in I}
V^*_{\o_{i^*}}\otimes\pi^*\CL_{-\o_i}^\phi)\otimes\omega_C)=
H^1(\Cx,(\bigoplus_{i\in I}
\CV^{*\CF}_{\o_{i^*}}\otimes\CL_{-\o_i}^\phi)\otimes\pi_*\omega_C)
\stackrel{d'_2}{\longrightarrow}$$
$$\stackrel{d'_2}{\longrightarrow}
H^0(\Cx,(\bigoplus_{i\in I}\CV^\CF_{\o_{i^*}}\otimes
\CL_{\o_i}^\phi)\otimes\pi_*\CO_C)=
H^0(C,\bigoplus_{i\in I}V_{\o_{i^*}}\otimes
\pi^*\CL_{\o_i}^\phi).$$
The bottom line of the present Section is the following

\cor{botline}
The following diagram commutes:

$\begin{CD}
H^1(C,(\bigoplus_{i\in I}
V^*_{\o_{i^*}}\otimes\pi^*\CL_{-\o_i}^\phi)\otimes\omega_C) @>>>
H^1(C,(\bigoplus_{i\in I}V^*_{\o_{i^*}}\otimes
\pi^*\CL_{-\o_i}^\phi)\otimes\Xi) \\
@V{d'_2}VV                @V{d^C_2}VV \\
H^0(C,\bigoplus_{i\in I}V_{\o_{i^*}}\otimes
\pi^*\CL_{\o_i}^\phi) @<<<
H^0(C,(\bigoplus_{i\in I}V_{\o_{i^*}}\otimes
\pi^*\CL_{\o_i}^\phi)\otimes\Xi) \\
\end{CD}$ \qed
\ecor

\ssec{c on C}{Calculation on $C$}
The differential $d^C_2$ of~\refc{botline} was computed in~\cite{fkmm}. To
formulate the result, we introduce homogeneous coordinates $z_1,z_2$ on $C$
such that $z=z_1/z_2$, so that $z_1=0$ (resp. $z_2=0$) is an equation of
$0\in C$ (resp. $\infty\in C$). We also introduce another copy of the curve
$C$ with homogeneous coordinates $u_1,u_2$, and $u:=u_1/u_2$. The differential
$d^C_2$ has ``matrix elements" $$\widetilde{D}_{\o_i,\o_j}\colon
H^1(C,V^*_{\o_{i^*}}\otimes\pi^*\CL_{-\o_i}^\phi\otimes\Xi)\to
H^0(C,V_{\o_{j^*}}\otimes\pi^*\CL_{\o_j}^\phi\otimes\Xi).$$ Note that
$\pi^*\CL_{-\o_i}\cong\CO_C(-a_i)$, and
$\pi^*\CL_{\o_j}\cong\CO_C(a_j)$, while $\Xi\cong\CO_C(-1)\oplus\CO_C(-1)$
is Serre selfdual, so that $$\widetilde{D}_{\o_i,\o_j}\in
V_{\o_{i^*}}\otimes V_{\o_{j^*}}\otimes
H^0(C,\CO_C(a_i-1)\oplus\CO_C(a_i-1))\otimes
H^0(C,\CO_C(a_j-1)\oplus\CO_C(a_j-1)).$$ Decomposing
$V_{\o_{i^*}}\otimes V_{\o_{j^*}}$ according to the weights of $T$, for
$\check\lambda\in\Lambda^\vee$ we obtain a matrix element
$\widetilde{D}_{\o_i,\o_j}^{\check\lambda}$ which is defined as the weight
$\check\lambda^*$-component of $\widetilde{D}_{\o_i,\o_j}$. Then according
to~\cite[3.8]{fkmm},
$\widetilde{D}_{\o_i,\o_j}^{\o_i+\o_j}=0=
\widetilde{D}_{\o_i,\o_i}^{2\o_i-2\check\alpha_i}$, and if $i\ne j$, then
$\widetilde{D}_{\o_i,\o_j}^{\o_i+\o_j-\check\alpha_i}=0=
\widetilde{D}_{\o_i,\o_j}^{\o_i+\o_j-\check\alpha_j}$, while
\eq{7}
\widetilde{D}_{\o_i,\o_i}^{2\o_i-\check\alpha_i}=\check{d}_i
(F_{i^*}v_{\o_{i^*}}\otimes v_{\o_{i^*}}-v_{\o_{i^*}}\otimes
F_{i^*}v_{\o_{i^*}})\otimes\frac{Q_i(z_1,z_2)R_i(u_1,u_2)-R_i(z_1,z_2)Q_i(u_1,u_2)}
{z_1u_2-z_2u_1}
\end{equation}
\begin{multline}\label{8}
\widetilde{D}_{\o_i,\o_j}^{\o_i+\o_j-\check\alpha_i-\check\alpha_j}=\\
(\check\alpha_i,\check\alpha_j)(\langle\alpha_i,\check\alpha_j\rangle^{-1}
v_{\o_{i^*}}\otimes F_{i^*}F_{j^*}v_{\o_{j^*}}+F_{i^*}v_{\o_{i^*}}\otimes
F_{j^*}v_{\o_{j^*}}+\langle\alpha_j,\check\alpha_i\rangle^{-1}
F_{j^*}F_{i^*}v_{\o_{i^*}}\otimes v_{\o_{j^*}})\otimes\\
\frac{R_i(z_1,z_2)R_j(u_1,u_2)-Q_i(z_1,z_2)S_{ji}(u_1,u_2)-
S_{ij}(z_1,z_2)Q_j(u_1,u_2)}{z_1u_2-z_2u_1}
\end{multline}
Here the homogeneous polynomials $Q_i,R_i,S_{ij}$ are but the homogenizations
of the same named polynomials of one variable introduced in~\refss{coord}, and
the above matrix coefficients are ``scalar" $2\times2$-matrices with respect
to the decomposition $\Xi\cong\CO_C(-1)\oplus\CO_C(-1)$.

Now to compute the desired $d'_2$ it remains to describe the horizontal arrows
of the commutative diagram of~\refc{botline}. The lower one
$$H^0(C,(\bigoplus_{i\in I}V_{\o_{i^*}}\otimes\pi^*\CL_{\o_i}^\phi)\otimes\Xi)\to
H^0(C,\bigoplus_{i\in I}V_{\o_{i^*}}\otimes\pi^*\CL_{\o_i}^\phi)$$
arises from the surjection $$H^0(C,\CO_C(a_i-1)\oplus\CO_C(a_i-1))\cong
H^0(C,\pi^*\CL^\phi_{\o_i}\otimes\Xi)\to H^0(C,\pi^*\CL^\phi_{\o_i})\cong
H^0(C,\CO_C(a_i))$$ which takes a pair $(P_1(z_1,z_2),P_2(z_1,z_2))$ of
homogeneous degree $a_i-1$ polynomials to $z_1P_1+z_2P_2$. The upper arrow
$$H^1(C,(\bigoplus_{i\in I}
V^*_{\o_{i^*}}\otimes\pi^*\CL_{-\o_i}^\phi)\otimes\omega_C)\to
H^1(C,(\bigoplus_{i\in I}V^*_{\o_{i^*}}\otimes
\pi^*\CL_{-\o_i}^\phi)\otimes\Xi)$$ arises from the dual embedding
$$H^0(C,\CO_C(a_i))^*\cong H^1(C,\pi^*\CL_{-\o_i}^\phi\otimes\omega_C)\to
H^1(C,\pi^*\CL_{-\o_i}^\phi\otimes\Xi)\cong
H^0(C,\CO_C(a_i-1)\oplus\CO_C(a_i-1))^*.$$ Namely, if we think of
$H^0(C,\CO_C(a_i))^*$ as of the homogeneous degree $a_i$ differential operators
in $(u_1,u_2)$, then an operator $O$ goes to $(O_1,O_2)$ such that
$O_1(P):=O(u_2\cdot P)/2$, while $O_2(P):=O(u_1\cdot P)/2$.

Composing with the matrix elements of~\refe{7},~(\ref{8}) we obtain the
corresponding matrix elements $'\!\widetilde{D}_{\o_i,\o_j}\in
V_{\o_{i^*}}\otimes V_{\o_{j^*}}\otimes
H^0(C,\CO_C(a_i))\otimes H^0(C,\CO_C(a_j))$ of $d'_2$:
\begin{multline}\label{77}
'\!\widetilde{D}_{\o_i,\o_i}^{2\o_i-\check\alpha_i}=\check{d}_i
(F_{i^*}v_{\o_{i^*}}\otimes v_{\o_{i^*}}-v_{\o_{i^*}}\otimes
F_{i^*}v_{\o_{i^*}})\otimes\\
\frac{(z_1u_2+z_2u_1)(Q_i(z_1,z_2)R_i(u_1,u_2)-R_i(z_1,z_2)Q_i(u_1,u_2))}
{2(z_1u_2-z_2u_1)}
\end{multline}
\begin{multline}\label{88}
'\!\widetilde{D}_{\o_i,\o_j}^{\o_i+\o_j-\check\alpha_i-\check\alpha_j}=\\
(\check\alpha_i,\check\alpha_j)(\langle\alpha_i,\check\alpha_j\rangle^{-1}
v_{\o_{i^*}}\otimes F_{i^*}F_{j^*}v_{\o_{j^*}}+F_{i^*}v_{\o_{i^*}}\otimes
F_{j^*}v_{\o_{j^*}}+\langle\alpha_j,\check\alpha_i\rangle^{-1}
F_{j^*}F_{i^*}v_{\o_{i^*}}\otimes v_{\o_{j^*}})\otimes\\
\frac{(z_1u_2+z_2u_1)(R_i(z_1,z_2)R_j(u_1,u_2)-Q_i(z_1,z_2)S_{ji}(u_1,u_2)-
S_{ij}(z_1,z_2)Q_j(u_1,u_2))}{2(z_1u_2-z_2u_1)}
\end{multline}

Going back from the homogeneous polynomials in $(z_1,z_2)$ (resp. $(u_1,u_2)$)
to the polynomials in $z=z_1/z_2$ (resp. $u=u_1/u_2$) we arrive at the
following

\prop{atlast}
The matrix elements of the differential $d'_2$~\refe{d'} are
\eq{777}
'\!\widetilde{D}_{\o_i,\o_i}^{2\o_i-\check\alpha_i}=\check{d}_i
(F_{i^*}v_{\o_{i^*}}\otimes v_{\o_{i^*}}-v_{\o_{i^*}}\otimes
F_{i^*}v_{\o_{i^*}})\otimes
\frac{(z+u)(Q_i(z)R_i(u)-R_i(z)Q_i(u))}{2(z-u)}
\end{equation}
\begin{multline}\label{888}
'\!\widetilde{D}_{\o_i,\o_j}^{\o_i+\o_j-\check\alpha_i-\check\alpha_j}=\\
(\check\alpha_i,\check\alpha_j)(\langle\alpha_i,\check\alpha_j\rangle^{-1}
v_{\o_{i^*}}\otimes F_{i^*}F_{j^*}v_{\o_{j^*}}+F_{i^*}v_{\o_{i^*}}\otimes
F_{j^*}v_{\o_{j^*}}+\langle\alpha_j,\check\alpha_i\rangle^{-1}
F_{j^*}F_{i^*}v_{\o_{i^*}}\otimes v_{\o_{j^*}})\otimes\\
\frac{(z+u)(R_i(z)R_j(u)-Q_i(z)S_{ji}(u)-
S_{ij}(z)Q_j(u))}{2(z-u)}
\end{multline}
while $'\!\widetilde{D}_{\o_i,\o_j}^{\o_i+\o_j}=0=\,
'\!\widetilde{D}_{\o_i,\o_i}^{2\o_i-2\check\alpha_i}$, and if $i\ne j$, then
$'\!\widetilde{D}_{\o_i,\o_j}^{\o_i+\o_j-\check\alpha_i}=0=\,
'\!\widetilde{D}_{\o_i,\o_j}^{\o_i+\o_j-\check\alpha_j}$.
\eprop

\ssec{calc}{Calculation of the Poisson bracket}

The differential $d_2$~(\refss{bivec}) defines a bivector field on $\oY^{\alp}$ (i.e. a bidifferential operation on the coordinate ring of $\oY^{\alp}$). We denote the bivector by $\fB$ and the corresponding bidifferential operation on the coordinate ring of $\oY^\alpha$ simply by $\{\cdot,\cdot\}$ (though it is \emph{not} a Poisson bracket on $\oY^\alp$).

\prop{bracket} We have $$\{w_{i,r},w_{j,s}\}=0,$$
$$ \{w_{i,r},y_{j,s}\}=\check d_i\delta_{ij}\delta_{rs}w_{i,r}y_{j,s},
$$
$$ \{y_{i,r},y_{j,s}\}=(1-\delta_{ij})(\check\alpha_i,\check\alpha_j)\frac{w_{i,r}+w_{j,s}}{2(w_{i,r}-w_{j,s})}y_{i,r}y_{j,s}.$$
\eprop

\prf
By~\refp{atlast} on $\oY^{\alp}$ we have
\begin{equation}\label{QQ}\{Q_i(z),Q_j(u)\}=0,\end{equation}
\begin{equation}\label{QR} \{Q_i(z),R_j(u)\}=
-\check{d}_i\delta_{ij}\frac{z+u}{2(z-u)}(Q_i(z)R_j(u)-
R_i(z)Q_j(u)),
\end{equation}
\begin{multline}\label{RR} \{R_i(z),R_j(u)\}=\\=(1-\delta_{ij})(\check\alpha_i,\check\alpha_j)
\frac{z+u}{2(z-u)}(R_i(z)R_j(u)- Q_i(z)S_{ji}(u)-S_{ij}(z)Q_j(u)).
\end{multline}
The relation $\{w_{i,r},w_{j,s}\}=0$ is obvious from~(\ref{QQ}). We have $\{B_i,B_j\}=0$ and $\{w_{i,r},B_j\}=0$ from~(\ref{QQ}) as well. Substituting $u=w_{j,s}$ to~(\ref{QR}), we get $\{B_i,y_{j,s}\}=-\frac{\check d_i\delta_{ij}}{2}B_iy_{j,s}$ and $\{w_{i,r},y_{j,s}\}=\check d_i\delta_{ij}\delta_{rs}w_{i,r}y_{j,s}$.
Finally, substituting $z=w_{i,r}, u=w_{j,s}$ to~(\ref{RR}), we get $\{y_{i,r},y_{j,s}\}=(1-\delta_{ij})(\check\alpha_i,\check\alpha_j)
\frac{w_{i,r}+w_{j,s}}{2(w_{i,r}-w_{j,s})}y_{i,r}y_{j,s}$.
\epr

The $1\times T$ action on $\oY^{\alp}$ preserves this bivector field, hence it gives a well-defined bivector field on $\oY^{\alp}/(1\times T)$. Moreover, the following is true:

\cor{quotient} The map $\varpi:\oY^{\alp}\to\dZ^\alp$ agrees with the bivector field $\fB$ on $\oY^\alp$ (in the sense that for $f_1,f_2\in\BC[\dZ^\alp]$ we have $\{\varpi^*(f_1),\varpi^*(f_2)\}=\varpi^*(f)$ for some $f\in\BC[\dZ^\alp]$). So we get a bivector field on $\dZ^\alpha=\varpi(\oY^{\alp})$.
\ecor

\prf
Note that the functions $w_{i,r},y_{j,s}$ form a (rational \'etale) coordinate system on $\dZ^\alpha$. So the only thing to be checked is that the bracket of any pair of pullbacks of these functions is a pullback of some function on $\dZ^\alpha$. But this immediately follows from \refp{bracket}.
\epr

Slightly abusing notation we denote the image of $\fB$ on $\dZ^\alp$ also by $\fB$.

\cor{bivector-Z} In the coordinates $w_{i,r}, y_{j,s}$ on $\dZ^\alpha$ the bivector field reads  $$\fB=\sum_{i,r}\check{d}_iw_{i,r}y_{i,r}\frac{\partial}{\partial w_{i,r}}\wedge \frac{\partial}{\partial y_{i,r}}+
\sum_{i\ne j}\sum_{r,s}\frac{(\check\alpha_i,\check\alpha_j)}
{2}
\frac{w_{i,r}+w_{j,s}}{w_{i,r}-w_{j,s}}y_{i,r}y_{j,s}\frac{\partial}{\partial y_{i,r}}\wedge \frac{\partial}{\partial y_{j,s}}.$$
\ecor

\cor{poisson} The bivector field $\fB$ on $\dZ^\alpha$ is Poisson, i.e.
$[\fB,\fB]=0$. This Poisson structure extends uniquely to $\hZ^\alpha$.
\ecor

\prf
The first claim is immediate from the explicit formula of~\refc{bivector-Z}.
We have a smooth open subvariety
$\dZ^\alpha\subset\buZ^\alpha\subset\hZ^\alpha$ formed by the
based quasimaps with defect at most a simple coroot, see
e.g.~\cite[proof of~Proposition~5.1]{bf}. Its complement has codimension 2 in
$\hZ^\alpha$. Now $\hZ^\alpha$ is normal by~\cite[Corollary~2.10]{bf}, so it
suffices to check that the bivector field on $\dZ^\alpha$
extends as a Poisson structure to $\buZ^\alpha$. Moreover,
it suffices to check this at the generic points of the boundary components
$\buZ^\alpha\setminus\dZ^\alpha$
(given by equations $y_{i,r}=0$) where the claim is evident
from the explicit formula of~\refc{bivector-Z}.
\epr

\cor{symplectic} The Poisson structure $\fB$ on $\dZ^\alpha$ is nondegenerate. The corresponding symplectic form reads  $$\Omega_\trig:=\fB^{-1}=\sum_{i,r}\frac{dy_{i,r}\wedge dw_{i,r}}{\check{d}_iw_{i,r}y_{i,r}}+
\sum_{i\ne j}\sum_{r,s}\frac{(\check\alpha_i,\check\alpha_j)}
{2\check{d}_i\check{d}_j}
\frac{w_{i,r}+w_{j,s}}{w_{i,r}-w_{j,s}}\frac{dw_{i,r}\wedge
dw_{j,s}}{w_{i,r}w_{j,s}}.$$
\ecor

\sec{slice}{Transversal slices in the affine flag variety}

\ssec{Waff}{Schubert cells in the affine flag varieties}
We have an embedding of the affine Grassmannian of $G$ (thin one: an ind-scheme)
into the Kashiwara affine Grassmannian (thick one: an infinite type scheme):
$\Gr=G((z))/G[[z]]\cong G[z^{\pm1}]/G[z]\hookrightarrow G((z^{-1}))/G[z]=\bGr$.
The subgroup of currents $G[[z]]$ (resp. $G[[z^{-1}]]$) taking value in $B$
(resp. $B_-$) at $z=0$ (resp. $z=\infty$) is the Iwahori group $\bIw$ (resp.
$\bIw_-$). The unipotent radical of $\bIw$ (resp. $\bIw_-$) is denoted
$\bN$ (resp. $\bN_-$). We have an embedding of the affine flag variety of $G$
(thin one: an ind-scheme)
into the Kashiwara affine flag variety (thick one: an infinite type scheme):
$\Fl=G((z))/\bIw\cong G[z^{\pm1}]/(G[z]\cap\bIw)\hookrightarrow
G((z^{-1}))/(G[z]\cap\bIw)=\bFl$. The natural projection
$\pr\colon \bFl\to\bGr$ (as well as its restriction $\pr\colon \Fl\to\Gr$) is a
fibration with fibers $\CB$.

The set of $T$-fixed points in $\Fl$ (resp. $\Gr$) is in a natural bijection
with the affine Weyl group $W_a=W\ltimes\Lambda$ (resp. the coweight lattice
$\Lambda$). For $w\in W_a$ we will denote the corresponding $T$-fixed point by
the same symbol $w$; its $\bN$-orbit (resp. $\bN_-$-orbit) will be denoted by
$\Fl^w\subset\Fl$ (resp. $\bFl_w\subset\bFl$): a thin (resp. thick) Schubert
cell. The intersection $\Fl^w_y:=\Fl^w\cap\bFl_y$ (an open Richardson variety,
aka transversal slice) is nonempty iff $w\geq y$ in the Bruhat order.
Similarly, for $\lambda\in\Lambda\subset\Gr$ the $\bN$-orbit $\bN\cdot\lambda$
(resp. $\bN_-$-orbit) will be denoted by $X^\lambda\subset\Gr$
(resp. $\bX_\lambda\subset\bGr$): a thin (resp. thick) Schubert cell.
For a dominant coweight $\lambda\in\Lambda^+\subset\Lambda$ the $G((z))$-orbit
$\Gr^\lambda:=G((z))\cdot\lambda$ is a union
$\Gr^\lambda=\bigsqcup_{\nu\in W\cdot\lambda}X^\nu$.

Recall the notations of~\cite[2.4]{bf}: $G_1\subset G[[z^{-1}]]$ is the kernel
of evaluation at $z=\infty$, and $\CW_\mu:=G_1\cdot\mu\subset\bGr$ for
$\mu\in\Lambda$. If $\mu$ is dominant, then $\CW_\mu=\bX_\mu$.
If $\lambda\geq\mu$ is also dominant, then the transversal slice
$\CW^\lambda_\mu:=\Gr^\lambda\cap\CW_\mu$ of {\em loc. cit.} is a union
$\CW^\lambda_\mu=\bigsqcup_{\nu\in W\cdot\lambda}X^\nu\cap\bX_\mu$.

Given a dominant $\eta\in\Lambda^+\subset\Lambda\subset W_a$, we consider
$-\eta$ as an element
in the affine Weyl group; it is the {\em minimal} length representative of
its {\em left} $W$-coset, and the {\em maximal} length representative of
its {\em right} $W$-coset. Furthermore, $\eta$ is the maximal length
representative of its left $W$-coset, and the minimal length representative
of its right $W$-coset.
The projection $\pr\colon W_a=\Fl^T\to\Gr^T=\Lambda$ realizes $\Lambda$ as the set
of {\em left} $W$-cosets in $W_a$. Hence for a dominant $\lambda\in W_a/W$,
the affine Weyl group element $w_0(\lambda)\in\Lambda\subset W_a$
(resp. $w_0(\lambda)\times w_0=w_0\times\lambda\in W\ltimes\Lambda=W_a$)
is the minimal length
(resp. maximal length) representative of the left $W$-coset $\lambda$.
In particular, $\pr\colon \Fl^{w_0(\lambda)}\iso X^\lambda$, and
$\pr\colon \bFl_{w_0\times\lambda}\iso\bX_\lambda$.

Finally, for $\lambda\geq\mu\in\Lambda^+$, the intersection
$X^\lambda\cap\bX_\mu$ is open in the slice $\CW^\lambda_\mu\subset\Gr$,
and $\pr$ is an open embedding of the open Richardson variety
$\Fl^{w_0\times\lambda}_{w_0\times\mu}=\Fl^{w_0\times\lambda}\cap
\bFl_{w_0\times\mu}$ into $X^\lambda\cap\bX_\mu$. All in all,
\eq{hook}
\pr\colon \Fl^{w_0\times\lambda}_{w_0\times\mu}\hookrightarrow\CW^\lambda_\mu
\end{equation}
is an open embedding.

\ssec{modul}{A modular interpretation}
Recall the morphism $s^{\lambda^*}_{\mu^*}\colon \CW^{\lambda^*}_{\mu^*}\to Z^\alpha$
of~\cite[Theorem 2.8]{bf} (here $\alpha=\lambda-\mu,\ \lambda^*=-w_0\lambda,\
\mu^*=-w_0\mu$). Recall the open
subscheme of periodic monopoles $\dZ^\alpha\subset Z^\alpha$ introduced
in~\refd{heart}.

\prop{coinc} The composition
$s^{\lambda^*}_{\mu^*}\circ\pr\colon \Fl^{w_0\times\lambda^*}_{w_0\times\mu^*}\to Z^\alpha$
is an open embedding with the image $\dZ^\alpha\subset Z^\alpha$.
\eprop

\prf
Recall from the proof of~\cite[Theorem~2.8]{bf} that the slice closure
$\overline\CW{}^{\lambda^*}_{\mu^*}=\bigsqcup_{\mu^*\leq\nu^*\leq\lambda^*}\CW^{\nu^*}_{\mu^*}$
is the moduli space of the following data:
$(\CF_{\on{triv}}\stackrel{\sigma}{\longrightarrow}\CF_G)$ where
$\CF_G$ is a $G$-bundle of isomorphism class $\mu^*$, and $\sigma$ is
an isomorphism from the trivial $G$-bundle away from $0\in C$ with a pole
of degree $\lambda^*$ at 0, such that the value of the Harder-Narasimhan flag
of $\CF_G$ at $\infty\in C$ is compatible with the complete flag
$\sigma(B_-)$. The bundle $\CF_G$ has a unique complete flag
($B$-structure) $\phi$ of degree $w_0\mu^*=-\mu$ with value $B_-\in\CB$ at
$\infty\in C$ (with respect to the trivialization $\sigma$ at $\infty$).
This flag can be transformed via $\sigma^{-1}$ to obtain a degree $\alpha$
generalized $B$-structure $\sigma^{-1}\phi$ in $\CF_{\on{triv}}$ without a pole
but possibly with a defect at $0\in C$. The morphism
$s^{\lambda^*}_{\mu^*}\colon \overline\CW{}^{\lambda^*}_{\mu^*}\to Z^\alpha$ takes
$(\CF_{\on{triv}}\stackrel{\sigma}{\longrightarrow}\CF_G)$ to $\sigma^{-1}\phi$.
The open subset $U\subset\overline\CW{}^{\lambda^*}_{\mu^*}$ given by the
condition that $\sigma^{-1}\phi$ has no defect at $0\in C$, is mapped
isomorphically onto $\oZ^\alpha$. We have a still smaller open subset
$U'\subset U$ given by the condition that the fiber of $\sigma^{-1}\phi$
at $0\in C$ is transversal to the flag $B\in\CB$. The open subset $U'$
is mapped isomorphically onto $\dZ^\alpha$. Thus we have to check
$\pr\colon \Fl^{w_0\times\lambda^*}_{w_0\times\mu^*}\iso U'$.

Recall the semiinfinite orbit $\CS^{\lambda^*}$ (whose intersection with
$\CW^{\lambda^*}_{\mu^*}$ is dense in $\CW^{\lambda^*}_{\mu^*}$). It is formed by the data
$(\CF_{\on{triv}}\stackrel{\sigma}{\longrightarrow}\CF_G)$ such that the
transformation $\sigma\phi_{\on{triv}}$ of the trivial complete flag
with fibers $B\in\CB$ in $\CF_{\on{triv}}$ via $\sigma$ is a $B$-structure in
$\CF_G$ without defect at $0\in C$. Note that $U'$ lies inside
$\CS^{\lambda^*}\cap\overline\CW{}^{\lambda^*}_{\mu^*}$ and is given there by the
condition that the fibers of $\sigma\phi_{\on{triv}}$ and of $\phi$ at
$0\in C$ are transversal.
According to~\cite[Theorem~3.2,~(3.6)]{mv}, for $\nu^*<\lambda^*$ we have
$\CS^{\lambda^*}\cap\Gr^{\nu^*}=\emptyset$, and
$\CS^{\lambda^*}\cap\Gr^{\lambda^*}=X^{\lambda^*}$. It follows that
$\CS^{\lambda^*}\cap\overline\CW{}^{\lambda^*}_{\mu^*}=X^{\lambda^*}\cap\bX_{\mu^*}$.
It remains to check that the open subset
$\pr(\Fl^{w_0\times\lambda^*}_{w_0\times\mu^*})\subset X^{\lambda^*}\cap\bX_{\mu^*}$
is nothing but $U'$.

To this end recall the modular interpretation of our slices.
%For the sake of exposition we assume that both
%$\lambda$ and $\mu$ are {\em regular} dominant coweights. The argument in the
%general case is similar but requires introducing more notations.
First of all,
$\bGr$ is the moduli space of $G$-bundles $\CF_G$ on $C$ equipped with a
trivialization $\varsigma$ in the formal neighbourhood of $\infty\in C$. Second,
$\bFl$ is the moduli space of triples $(\CF_G,\varsigma,F)$ where
$(\CF_G,\varsigma)\in\bGr$, and $F$ is a $B$-structure in the fiber of $\CF_G$
at $0\in C$. Third, $\bX_{\mu^*}=\CW_{\mu^*}\subset\bGr$ is formed by the pairs
$(\CF_G,\varsigma)$
such that the isomorphism type of $\CF_G$ is $\mu^*$, and the value of the
Harder-Narasimhan flag of $\CF_G$ at $\infty\in C$ is compatible with
$B_-\in\CB$ (with respect to the
trivialization $\varsigma$ at $\infty\in C$).
Now $\bFl_{w_0\times\mu^*}\subset\bFl$
is formed by the triples $(\CF_G,\varsigma,F)$ such that
$(\CF_G,\varsigma)\in\bX_{\mu^*}$, and $F$ is the value $\phi|_0$ at $0\in C$
of the unique degree $w_0\mu^*=-\mu$ complete flag $\phi$
in $\CF_G$ such that $\phi|_\infty=B_-\in\CB$ (so that $\phi$ is the
refinement of the Harder-Narasimhan flag of $\CF_G$). Furthermore, $\Gr$ is the moduli space of
$G$-bundles $\CF_G$ on $C$ equipped with a trivialization $\sigma$ over
$C\setminus0$, while $\Fl$ is the moduli space of triples $(\CF_G,\sigma,F)$
where $(\CF_G,\sigma)\in\Gr$, and $F$ is a $B$-structure in the fiber of
$\CF_G$ at $0\in C$. The projection $\pr\colon \Fl\to\Gr$ admits a section
$s$ over $X^{\lambda^*}=\CS^{\lambda^*}\cap\Gr^{\lambda^*}$: we define $F$ as the fiber
at $0\in C$ of the transformation $\sigma\phi_{\on{triv}}$
of the trivial $B$-structure $B\in\CB$ in the trivial $G$-bundle.
Finally, $\Fl^{w_0\times\lambda^*}\subset\Fl$ is formed by the triples
$(\CF_G,\sigma,F)$ such that $(\CF_G,\sigma)\in X^{\lambda^*}$, and
$F$ is transversal to $s(\CF_G,\sigma)$.

Thus $\pr(\Fl^{w_0\times\lambda^*}_{w_0\times\mu^*})=
U'\subset X^{\lambda^*}\cap\bX_{\mu^*}$.
The proposition is proved.
\epr

\ssec{stabi}{Stabilization}
Let $\mu,\nu\in\Lambda^+$ be dominant coweights. According to~\cite[2E]{kwy},
we have the inclusion of stabilizers
$\on{St}_\mu\subset\on{St}_{\mu+\nu}\subset G_1\subset G[[z^{-1}]]$,
so the identity morphism $G_1\to G_1$ induces a morphism
$\varsigma_\mu^{\mu+\nu}\colon \bX_\mu=G_1/\on{St}_\mu\to G_1/\on{St}_{\mu+\nu}=
\bX_{\mu+\nu}$. According to {\em loc. cit.}, $\varsigma_\mu^{\mu+\nu}$ restricts
to the same named morphism $\CW^\lambda_\mu\to\CW^{\lambda+\nu}_{\mu+\nu}$ for any
$\Lambda^+\ni\lambda\geq\mu$.
Similarly, we have $\on{St}_{w_0\times\mu}\subset\on{St}_{w_0\times(\mu+\nu)}\subset
\bN_-$, and the identity morphism $\bN_-\to\bN_-$ induces a morphism
$\sigma_\mu^{\mu+\nu}\colon \bFl_{w_0\times\mu}=\bN_-/\on{St}_{w_0\times\mu}\to
\bN_-/\on{St}_{w_0\times(\mu+\nu)}=\bFl_{w_0\times(\mu+\nu)}$ which restricts to
the same named morphism $\Fl^{w_0\times\lambda}_{w_0\times\mu}\to
\Fl^{w_0\times(\lambda+\nu)}_{w_0\times(\mu+\nu)}$ for any
$\Lambda^+\ni\lambda\geq\mu$. The following diagram commutes:
\eq{squa}
\begin{CD}
\Fl^{w_0\times\lambda}_{w_0\times\mu} @>{\sigma_\mu^{\mu+\nu}}>>
\Fl^{w_0\times(\lambda+\nu)}_{w_0\times(\mu+\nu)}\\
@V{\pr}VV        @V{\pr}VV\\
\CW^\lambda_\mu @>{\varsigma_\mu^{\mu+\nu}}>>        \CW^{\lambda+\nu}_{\mu+\nu}\\
\end{CD}
\end{equation}
Moreover, from the construction of
$s^{\lambda^*}_{\mu^*}\colon \CW^{\lambda^*}_{\mu^*}\to Z^\alpha$
in~\cite[Lemma~2.7,~Theorem~2.8]{bf} (where $\alpha=\lambda-\mu$)
it follows immediately that the following diagrams commute as well:
\eq{squar}
\begin{CD}
\CW^{\lambda^*}_{\mu^*} @>{\varsigma_{\mu^*}^{\mu^*+\nu^*}}>>
\CW^{\lambda^*+\nu^*}_{\mu^*+\nu^*}\\
@V{s^{\lambda^*}_{\mu^*}}VV        @V{s^{\lambda^*+\nu^*}_{\mu^*+\nu^*}}VV\\
Z^\alpha @=        Z^\alpha\\
\end{CD}
\end{equation}

\eq{square}
\begin{CD}
\Fl^{w_0\times\lambda^*}_{w_0\times\mu^*} @>{\sigma_{\mu^*}^{\mu^*+\nu^*}}>>
\Fl^{w_0\times(\lambda^*+\nu^*)}_{w_0\times(\mu^*+\nu^*)}\\
@V{s^{\lambda^*}_{\mu^*}\circ\pr}VV      @V{s^{\lambda^*+\nu^*}_{\mu^*+\nu^*}\circ\pr}VV\\
\dZ^\alpha @=        \dZ^\alpha\\
\end{CD}
\end{equation}
It follows in particular that $\sigma_{\mu^*}^{\mu^*+\nu^*}\colon
\Fl^{w_0\times\lambda^*}_{w_0\times\mu^*}\iso
\Fl^{w_0\times(\lambda^*+\nu^*)}_{w_0\times(\mu^*+\nu^*)}$ is an isomorphism.

\ssec{morp}{$s^\lambda_\mu$ in coordinates}
We will use the generalized minors of~\cite[2A]{kwy} to construct regular
functions on the open Richardson varieties.
%however, our notation differs slightly from the one of {\em loc. cit.}
Namely, given an irreducible
$G$-module $V_{\check\lambda}$ with highest weight
$\check\lambda\in\Lambda^\vee_+$ and highest vector $v_{\check\lambda}$,
its dual $V_{\check\lambda}^*$ is isomorphic to $V_{\check\lambda^*}$
with the lowest weight $-\check\lambda$, and the lowest vector
$v_{-\check\lambda}$ such that
$\langle v_{-\check\lambda}, v_{\check\lambda}\rangle=1$.
Given $w,y\in W$, we define the following regular function on $G$
\eq{vardel}
\varDelta_{w\check\lambda,y\check\lambda}(g):=\langle\ol{w}v_{-\check\lambda},
g\ol{y}v_{\check\lambda}\rangle
\end{equation}
where $\ol{w},\ol{y}\in G$ are the lifts of $w,y$ defined in {\em loc. cit.}
%(note the difference with the generalized minor
%$\Delta_{w\check\lambda,y\check\lambda}(g)$ of {\em loc. cit.})

Following {\em loc. cit.} we consider the regular functions
$\varDelta_{w\check\lambda,y\check\lambda}^{(s)},\ s\in\BZ$, on $G((z^{-1}))$
defined as follows:
\eq{vardels}
\varDelta_{w\check\lambda,y\check\lambda}(g(z))=\sum_{s=-\infty}^{\infty}
\varDelta_{w\check\lambda,y\check\lambda}^{(s)}(g(z))z^{-s}
\end{equation}

More generally, to any $v\in V_{\check\lambda}$ and $\beta\in V_{\check\lambda}^*$ we can assign the generalized minor $\varDelta_{\beta,v}(z):=\langle\beta,
g(z)v\rangle$. We also denote by $\varDelta_{\beta,v}^{(s)}$ the coefficient at $z^{-s}$ of the power series $\varDelta_{\beta,v}(z)$.

Recall from~\cite[2.6]{bf} that $\on{St}_\mu\subset G_1\subset G[[z^{-1}]]$
is the stabilizer of $\mu\in\Lambda^+\subset\Lambda=\Gr^T$. Similarly,
for $w\in W_a$ we denote by $\on{St}_w\subset\bN_-$ the stabilizer of
$w\in W_a=\Fl^T$. We have $\bX_\mu=G_1/\on{St}_\mu$, and
$\bFl_w=\bN_-/\on{St}_w$. In case $w=w_0\times\mu$ (see~\refss{Waff}),
we have $\on{St}_\mu=G_1\cap\on{St}_{w_0\times\mu}$, and the natural morphism
$\bFl_{w_0\times\mu}=\bN_-/\on{St}_{w_0\times\mu}\to G_1/\on{St}_\mu=\bX_\mu$
is an isomorphism. According to~\cite[Lemma~2.19]{kwy}, the functions
$\varDelta_{\check\omega_i,\check\omega_i}^{(s)},\ \varDelta_{s_i\check\omega_i,\check\omega_i}^{(s)},\
s>0,\ i\in I$ restricted to $G_1$ (resp. $\bN_-$) are $\on{St}_\mu$-invariant
(resp. $\on{St}_{w_0\times\mu}$-invariant); hence they may be viewed as the
functions on $\bFl_{w_0\times\mu}\cong\bX_\mu$.

Now let $\Lambda^+\ni\lambda\geq\mu$, and
$\Lambda_+\ni\alpha=\sum_{i\in I}a_i\alpha_i:=\lambda-\mu$. Recall the
isomorphism $s^{\lambda^*}_{\mu^*}\circ\pr\colon \Fl^{w_0\times\lambda^*}_{w_0\times\mu^*}\iso
\dZ^\alpha$ of~\refp{coinc}. Recall also the regular polynomial-valued functions
$Q_i,R_i$ on $Z^\alpha$ (see e.g.~\cite[3.3]{fkmm}):
$Q_i=z^{a_i}+q_{i,a_i-1}z^{a_i-1}+\ldots$
(resp. $R_i= r_{i,a_i-1}z^{a_i-1}+\ldots$) is the highest (resp. prehighest)
Pl\"ucker coordinate on the space of based quasimaps (in notations of
{\em loc. cit.} $Q_i=\phi_{\o_i}^{-\o_i},\
R_i=\phi_{\o_i}^{\check\alpha_i-\o_i}$). Following {\em loc. cit.}
and~\refss{coord} we also consider a rational \'etale
coordinate system on $Z^\alpha$. Namely,
$(w_{i,r})_{i\in I}^{1\leq r\leq a_i}$ are the ordered roots of $Q_i$, and
$y_{i,r}:=R_i(w_{i,r})$.

\prop{comp}
Under the isomorphism $s^{\lambda^*}_{\mu^*}\circ\pr\colon
\Fl^{w_0\times\lambda^*}_{w_0\times\mu^*}\iso\dZ^\alpha$ we have
$Q_i=\sum_{s=0}^{a_i}\varDelta_{\check\omega_i,\check\omega_i}^{(s)}z^{a_i-s},\
R_i=\sum_{s=0}^{a_i}\varDelta_{s_i\check\omega_i,\check\omega_i}^{(s)}z^{a_i-s}$.
\eprop

\prf
Follows at once from the commutative diagram~\cite[(2.3)]{bf}
(and the definition of $\pi_{\mu,n}$ in~\cite[Lemma~2.7]{bf}).
\epr

\ssec{ratpois}{Rational Poisson bracket revisited}

We fix a basis $e_\alpha, e_{-\alpha}, h_i$ in $\fg$ where $i\in I$, and
$\alpha\in R^+$ is a positive coroot (and the weight of $e_\alpha$ is the dual
root $\check\alpha$; in particular, $e_{\alpha_i}=E_i$ of~\refss{BunT},
and $e_{-\alpha_i}=\check{d}_iF_i$).
We assume $(e_\alpha, e_{-\alpha})=1$, $(h_i, h_j)=\delta_{ij}$. Then the Lie bialgebra structure on $\fg((z^{-1}))$ is determined by the classical rational $r$-matrix
\eq{bd-rmatrix-rat}
r_\rat(z,u):=\frac{1}{z-u}(\sum\limits_{\alpha>0}e_\alpha\otimes e_{-\alpha}+e_{-\alpha}\otimes e_{\alpha}+\sum\limits_{i\in I}h_i\otimes h_i),
\end{equation}
see e.g.~\cite[Section~6.4]{es}. This determines a Poisson group structure on $G((z^{-1}))$ such that $G_1$ is a Poisson subgroup.

\prop{rational-RTT}\cite[Proposition~2.13]{kwy} The rational Poisson bracket
$\{,\}_\rat$ of the functions $\varDelta_{\beta,v}^{(s)}$ on the subgroup $G_1$ is
\begin{multline*}
\{\varDelta_{\beta_1,v_1}(z),\varDelta_{\beta_2,v_2}(u)\}_\rat=\\=
\frac{1}{z-u}(\sum\limits_{\alpha>0} \varDelta_{\beta_1,e_\alpha v_1}(z)\varDelta_{\beta_2,e_{-\alpha}v_2}(u)+\sum\limits_{\alpha>0} \varDelta_{\beta_1,e_{-\alpha} v_1}(z)\varDelta_{\beta_2,e_{\alpha}v_2}(u)+\sum\limits_{i\in I} \varDelta_{\beta_1, h_i v_1}(z)\varDelta_{\beta_2,h_iv_2}(u))-\\-
\frac{1}{z-u}(\sum\limits_{\alpha>0} \varDelta_{e_\alpha\beta_1, v_1}(z)\varDelta_{e_{-\alpha}\beta_2,v_2}(u)+\sum\limits_{\alpha>0} \varDelta_{e_{-\alpha}\beta_1, v_1}(z)\varDelta_{e_{\alpha}\beta_2,v_2}(u)+\sum\limits_{i\in I} \varDelta_{h_i \beta_1, v_1}(z)\varDelta_{h_i \beta_2,v_2}(u)).
\end{multline*}
\eprop

According to \cite{kwy}, this Poisson structure on $G_1$ induces a Poisson structure on transversal slices $\CW^\lambda_\mu$ in the affine Grassmannian $\bGr=G((z^{-1}))/G[z]$. On the other hand,
%for any transversal slice in
%$\bGr$ we have a map to the corresponding zastava space
%$s^\lambda_\mu\colon \CW^\lambda_\mu\to Z^\alpha$
%(where $\alpha=\lambda-\mu$) and a natural
recall a symplectic structure on $\oZ^\alpha$ defined in~\cite{fkmm}.
It extends uniquely to a Poisson bracket $\{,\}^Z_\rat$ on $Z^\alpha$ by the same
argument as in the proof of~\refc{poisson}.
%which extends as a Poisson structure to $Z^\alpha$ due to the following

%\lem{joel}
%There is a unique Poisson bracket $\{,\}^Z_\rat$ on
%$Z^\alpha$ whose restriction
%to $\oZ^\alpha\subset Z^\alpha$ comes from the symplectic structure
%of~\cite{fkmm}.
%\elem

%\prf
%We have a smooth open subvariety
%$\oZ^\alpha\subset\overset{\bullet}{Z}{}^\alpha\subset Z^\alpha$ formed by the
%based quasimaps with defect at most a simple coroot, see
%e.g.~\cite[proof of~Proposition~5.1]{bf}. Its complement has codimension 2 in
%$Z^\alpha$. Now $Z^\alpha$ is normal by~\cite[Corollary~2.10]{bf}, so it
%suffices to check that the symplectic structure of~\cite{fkmm} on $\oZ^\alpha$
%extends as a Poisson structure to $\overset{\bullet}{Z}{}^\alpha$. Moreover,
%it suffices to check this at the generic points of the boundary components
%$\overset{\bullet}{Z}{}^\alpha\setminus\oZ^\alpha$
%(given by equations $y_{i,r}=0$) where the claim is evident
%from the explicit formulas of~\cite[Proposition~2]{fkmm}.
%\epr

The following theorem confirms expectations of~\cite[Remark~2.11]{kwy}.

\th{rational-poisson-map} The map
$s^{\lambda^*}_{\mu^*}\colon \CW^{\lambda^*}_{\mu^*}\to Z^\alpha$ is Poisson.
\eth

\prf
The field of rational functions on $Z^\alpha$ coincides with the field of
rational functions in the Fourier coefficients of the functions $Q_i(z), R_i(z)$. Hence it is sufficient to show that the Poisson bracket of the coefficients $Q_i(z), R_i(z)$ is the same on $\CW^\lambda_\mu$ and $Z^\alpha$. Let us introduce the following generalized
minors:
$S_{ij}(z):=\varDelta_{E_jE_iv_{-\check\omega_i},v_{\check\omega_i}}=
\langle E_jE_iv_{-\check\omega_i},g(z)v_{\check\omega_i}\rangle$.
According to \cite[(7)~and~(8)]{fkmm}, the Poisson bracket of the (polynomial-valued)
functions $Q_i(z), R_i(z)$ is given by
\begin{equation}\label{rQQ}\{Q_i(z),Q_j(u)\}^Z_\rat=0,\end{equation}
\begin{equation}\label{rQR} \{Q_i(z),R_j(u)\}^Z_\rat=
-\check{d}_i\delta_{ij}(\frac{1}{z-u}Q_i(z)R_j(u)-
\frac{1}{z-u}R_i(z)Q_j(u)),
\end{equation}
\begin{multline}\label{rRR} \{R_i(z),R_j(u)\}^Z_\rat=\\=(1-\delta_{ij})((\check\alpha_i,\check\alpha_j)
\frac{1}{z-u}R_i(z)R_j(u)+\check d_i\check d_j\frac{1}{z-u}Q_i(z)S_{ji}(u))+\check d_i\check d_j\frac{1}{z-u}S_{ij}(z)Q_j(u)).
\end{multline}
On the other hand, the Fourier coefficients of the pullbacks
$(s^{\lambda^*}_{\mu^*})^* Q_i=
z^{a_i}+\sum_{s=1}^{a_i}\varDelta_{\check\omega_i,\check\omega_i}^{(s)}z^{a_i-s}$ and
$(s^{\lambda^*}_{\mu^*})^* R_i=
\sum_{s=1}^{a_i}\varDelta_{s_i\check\omega_i,\check\omega_i}^{(s)}z^{a_i-s}$
obey the same relations by \refp{rational-RTT}.
\epr

\ssec{pois}{Trigonometric Poisson bracket}

The standard Lie bialgebra structure on $\fg((z^{-1}))\oplus\ft$
(see e.g.~\cite[6.2.1,~6.5]{es}) gives rise to a Poisson structure on $\bFl$
such that the open Richardson varieties $\Fl^w_y$ are Poisson subvarieties
of $\bFl$ (cf.~\cite[Corollary~2.9]{ly}).

This Lie bialgebra structure on $\fg((z^{-1}))\oplus\ft$ is determined by the classical $r$-matrix
\eq{bd-rmatrix}
r_{\trig}(z,u):=\frac{1}{z-u}(z(\sum\limits_{\alpha>0}e_\alpha\otimes e_{-\alpha}+\frac{1}{2}\sum\limits_{i\in I}h_i\otimes h_i)+u(\sum\limits_{\alpha>0}e_{-\alpha}\otimes e_{\alpha}+\frac{1}{2}\sum\limits_{i\in I}h_i\otimes h_i)),
\end{equation}
see e.g.~\cite[(6.6)]{es}.

\prop{trig-RTT} The Poisson bracket of the functions $\varDelta_{\beta,v}^{(s)}$ on the Iwahori subgroup $\bIw_-$ is
\begin{multline*}
\{\varDelta_{\beta_1,v_1}(z),\varDelta_{\beta_2,v_2}(u)\}_{\trig}=\\=\frac{1}{z-u}(z\sum\limits_{\alpha>0} \varDelta_{\beta_1,e_\alpha v_1}(z)\varDelta_{\beta_2,e_{-\alpha}v_2}(u)+u\sum\limits_{\alpha>0} \varDelta_{\beta_1,e_{-\alpha} v_1}(z)\varDelta_{\beta_2,e_{\alpha}v_2}(u)+\\+\frac{z+u}{2}\sum\limits_{i\in I} \varDelta_{\beta_1, h_i v_1}(z)\varDelta_{\beta_2,h_iv_2}(u))-\\-
\frac{1}{z-u}(z\sum\limits_{\alpha>0} \varDelta_{e_\alpha\beta_1, v_1}(z)\varDelta_{e_{-\alpha}\beta_2,v_2}(u)+u\sum\limits_{\alpha>0} \varDelta_{e_{-\alpha}\beta_1, v_1}(z)\varDelta_{e_{\alpha}\beta_2,v_2}(u)+\\+\frac{z+u}{2}\sum\limits_{i\in I} \varDelta_{h_i \beta_1, v_1}(z)\varDelta_{h_i \beta_2,v_2}(u))
\end{multline*}
\eprop

\prf
This follows from the Belavin-Drinfeld formula for trigonometric $r$-matrix
(see e.g.~\cite[(6.6)]{es}).
Indeed, following \cite[Proposition~2.13]{kwy} we note that the cobracket on $\fg((z^{-1}))$ is coboundary, namely it is given by the map
$$
a(t)\mapsto[a(z)\otimes1+1\otimes a(u),r_{\trig}(z,u)],
$$
where the $r$-matrix is given by~\refe{bd-rmatrix}. By the standard procedure this gives a structure of Poisson group on $G((z^{-1}))$. We note that the
Iwahori subgroup $\bIw_-\subset G((z^{-1}))$ is a Poisson subgroup, hence the bracket of any two functions on it is the restriction of the bracket of any extensions of these functions to $G((z^{-1}))$. The rest of the proof is a word-to-word repetition of that of~\cite[Proposition~2.13]{kwy}.
\epr

By an abuse of notation, we will denote the
rational \'etale functions
$w_{i,r}\circ s^{\lambda^*}_{\mu^*}\circ\pr,\ y_{i,r}\circ s^{\lambda^*}_{\mu^*}\circ\pr$
(notations of~\refss{morp})
on $\Fl^{w_0\times\lambda^*}_{w_0\times\mu^*}$ simply by $w_{i,r},y_{i,r}$.

\prop{trig-xy} We have $$\{w_{i,r},w_{j,s}\}_{\trig}=0,$$
$$ \{w_{i,r},y_{j,s}\}_{\trig}=\check d_i\delta_{ij}\delta_{rs}w_{i,r}y_{j,s},
$$
$$ \{y_{i,r},y_{j,s}\}_{\trig}=(1-\delta_{ij})(\check\alpha_i,\check\alpha_j)\frac{w_{i,r}+w_{j,s}}{2(w_{i,r}-w_{j,s})}y_{i,r}y_{j,s}.$$
\eprop

\prf
Consider the functions $Q_i(z)=\varDelta_{\check\omega_i,\check\omega_i}^{(0)}\prod\limits_{r=1}^{a_i}(z-w_{i,r})$, $R_i(z)=\varDelta_{\check\omega_i,\check\omega_i}^{(0)}\sum\limits_{r=1}^{a_i}y_{i,r}\frac{Q_i(z)}{(z-w_{i,r})Q_i'(w_{i,r})}$.
According to \refp{comp} we have
$Q_i=\sum_{s=0}^{a_i}\varDelta_{\check\omega_i,\check\omega_i}^{(s)}z^{a_i-s},\
R_i=\sum_{s=0}^{a_i}\varDelta_{s_i\check\omega_i,\check\omega_i}^{(s)}z^{a_i-s}.$

Set $B_i:=\varDelta_{\check\omega_i,\check\omega_i}^{(0)}$ and recall
the generalized minors $S_{ij}(z)$ introduced in the proof of~\reft{rational-poisson-map}: $S_{ij}(z):=\varDelta_{E_jE_iv_{-\check\omega_i},v_{\check\omega_i}}=
\langle E_jE_iv_{-\check\omega_i},g(z)v_{\check\omega_i}\rangle$.
Then by~\refp{trig-RTT} we have
\begin{equation}\label{tQQ}\{Q_i(z),Q_j(u)\}_{\trig}=0,\end{equation}
\begin{equation}\label{tQR} \{Q_i(z),R_j(u)\}_{\trig}=
-\check{d}_i\delta_{ij}(\frac{z+u}{2(z-u)}Q_i(z)R_j(u)-
\frac{u}{z-u}R_i(z)Q_j(u)),
\end{equation}
\begin{multline}\label{tRR} \{R_i(z),R_j(u)\}_{\trig}=\\=(1-\delta_{ij})((\check\alpha_i,\check\alpha_j)\frac{z+u}{2(z-u)}R_i(z)R_j(u)+\check d_i\check d_j\frac{z}{z-u}Q_i(z)S_{ji}(u)+\check d_i\check d_j\frac{u}{z-u}S_{ij}(z)Q_j(u)).
\end{multline}
The relation $\{w_{i,r},w_{j,s}\}_{\trig}=0$ is obvious from~(\ref{tQQ}). Substituting $u=w_{j,s}$ to~(\ref{tQR}), we get $\{B_i,y_{j,s}\}_{\trig}=-\frac{\check d_i\delta_{ij}}{2}B_iy_{j,s}$ and $\{w_{i,r},y_{j,s}\}_{\trig}=\check d_i\delta_{ij}\delta_{rs}w_{i,r}y_{j,s}$.
Finally, substituting $z=w_{i,r}, u=w_{j,s}$ to~(\ref{tRR}), we get $\{y_{i,r},y_{j,s}\}_{\trig}=(1-\delta_{ij})(\check\alpha_i,\check\alpha_j)
\frac{w_{i,r}+w_{j,s}}{2(w_{i,r}-w_{j,s})}y_{i,r}y_{j,s}$.
\epr

\th{match}
The isomorphism $s^{\lambda^*}_{\mu^*}\circ\pr\colon
\Fl^{w_0\times\lambda^*}_{w_0\times\mu^*}\iso\dZ^\alpha$ of~\refp{coinc}
is a symplectomorphism.
\eth

\prf
Indeed, by \refp{trig-xy} and \refp{bracket} the Poisson brackets
$\{\cdot,\cdot\}$ on $\dZ^\alpha$ and $\{\cdot,\cdot\}_{\trig}$ on
$\Fl^{w_0\times\lambda^*}_{w_0\times\mu^*}$ are given by the same formulas on
coordinate functions $w_{i,r},\ y_{i,r}$.
\epr

\rem{ploxo} Note that the formulas (\ref{tQR}) and (\ref{tRR}) are different from (\ref{QR}) and (\ref{RR}), so the morphism $s^{\lambda^*}_{\mu^*}\circ\pr\colon
\Fl^{w_0\times\lambda^*}_{w_0\times\mu^*}\iso\dZ^\alpha$ \emph{does not} extend to a Poisson morphism $\bIw_-\to \oY^\alp$.
\erem

\sec{clus}{A speculation on cluster structure}

\ssec{gaff}{An affine Lie algebra}
Let $\hg$ be the universal central extension of the polynomial loop algebra
$\fg[z^{\pm1}]$:
\eq{centr}
0\to\BC\to\hg\to\fg[z^{\pm1}]\to0
\end{equation}
Let $\gaff=\hg\rtimes\BC d$ be the semidirect product of $\hg$ with the
degree operator. Then $\gaff$ is an untwisted affine Kac-Moody
Lie algebra. It has a triangular decomposition
$\gaff=\fn_-\oplus\ft_{\on{aff}}\oplus\fn$ where
$\fn_-=\on{Lie}\bN_-\cap\fg[z^{-1}],\
\fn=\on{Lie}\bN\cap\fg[z]$, and $\ft_{\on{aff}}$ is the affine Cartan
subalgebra. The fundamental weights will be denoted
$\varpi_i,\ i\in I_a:=I\sqcup\{i_0\}$. The corresponding fundamental
integrable representations (where $\fn$ acts locally nilpotently) will be
denoted $V_{\varpi_i}$, and their restricted duals (where $\fn_-$ acts
locally nilpotently) will be denoted $V_{\varpi_i}^*$. We choose the highest
weight vectors $v_{\varpi_i}\in V_{\varpi_i}$ and the lowest weight vectors
$v_{-\varpi_i}\in V_{\varpi_i}^*$ such that
$\langle v_{-\varpi_i},v_{\varpi_i}\rangle=1$.
Note that the action of $\fn_-$ (resp. $\fn$) on $V_{\varpi_i}^*$
(resp. $V_{\varpi_i}$) integrates to the action of $\bN_-$ (resp. $\bN$).
Given $w,y\in W_a$ and $i\in I_a$ we define the following regular function on
$\bN_-$ (a generalized minor):
\eq{vardelt}
\varDelta_{w\varpi_i,y\varpi_i}(g):=\langle\ol{w}v_{-\varpi_i},
g\ol{y}v_{\varpi_i}\rangle
\end{equation}
where $\ol{w},\ol{y}\in G_{\on{aff}}$ are the lifts of $w,y$ defined similarly
to~\cite[2A]{kwy}.

\ssec{seed}{An initial seed}
B.~Leclerc defines in~\cite{l} a cluster structure on the open Richardson
varieties in the flag varieties of simple Lie algebras of types $ADE$, but
presumably the construction can be extended to the affine Lie algebras of
arbitrary types. Here we describe the initial seed for
$\Fl_{w_0}^{w_0\times\lambda}$ following~\cite[Section~5,~4.8.3,~Corollary~4.4]{l}.
We choose a reduced expression in the affine Weyl group:
$\lambda=s_{i_1}\ldots s_{i_l}$ where $i_1,\ldots,i_l\in I_a$, and
$l=2|\lambda|$ is the length of
$\lambda\in\Lambda^+\subset\Lambda_+\subset\Lambda\subset W_a$
(for $\lambda=\sum_{i\in I}a_i\alpha_i$ we have $2|\lambda|=2\sum_{i\in I}a_i$).
Note that $i_1=i_0$ (the affine simple reflection).
Then the initial seed consists of all the (irreducible factors of the)
generalized minors
$\varDelta_{w_0s_{i_1}\cdots s_{i_r}\varpi_{i_r},w_0\varpi_{i_r}},\
1\leq r\leq l$ (they are well defined as functions on
$\Fl_{w_0}^{w_0\times\lambda}$ according to {\em loc. cit.}). Among them,
those which divide
$\prod_{i\in I_a}\varDelta_{(w_0\times\lambda)\varpi_i,w_0\varpi_i}$ are the
frozen variables.

\ssec{exch}{An exchange matrix}
The rows of the exchange matrix $B$ are numbered by $1\leq r\leq l$, and the
columns are numbered by those $1\leq s\leq l$ for which there exists
$r>s$ such that $i_r=i_s$ (the minimal among such $r$ is denoted $s^+$).
The matrix entries are as follows: $b_{s,s^+}=-b_{s^+,s}=-1$; and for
$s<r<s^+$ such that for any $r<r'<s^+$ we have $i_r\ne i_{r'}$, the matrix
entry $b_{s,r}=-C_{i_s,i_r}$, and $b_{r,s}=C_{i_r,i_s}$ (here
$(C_{i,j})_{i,j\in I_a}$ is the Cartan matrix of $\hg$). All the other
matrix entries are zero.

According to~\cite[Section~6]{l}, this cluster structure on
$\Fl_{w_0}^{w_0\times\lambda}$ is compatible with the symplectic structure
of~\refss{pois} on $\Fl_{w_0}^{w_0\times\lambda}$ in the sense
of~\cite[Section~4.1]{gsv}. Taking $\mu=0$ and $\alpha=\lambda$,
and transferring the cluster structure via
the isomorphism $s^{\lambda^*}_0\circ\pr\colon
\Fl^{w_0\times\lambda^*}_{w_0}\iso\dZ^\alpha$ we obtain a cluster
structure on $\dZ^\alpha$ compatible with its symplectic structure
(see~\reft{match}).

\ssec{desta}{Destabilization}
Let $\nu\in\Lambda^+$ be a dominant coweight. Then the open Richardson
variety $\Fl_{w_0\times\nu}^{w_0\times(\lambda+\nu)}$ also has a cluster
structure with the initial cluster given by certain generalized minors, and
with the same exchange matrix as in~\refss{exch}. However, the stabilization
map~\refe{square} $\sigma_0^\nu\colon \Fl_{w_0}^{w_0\times\lambda}\iso
\Fl_{w_0\times\nu}^{w_0\times(\lambda+\nu)}$ {\em does not} take the
initial seed of $\Fl_{w_0}^{w_0\times\lambda}$ to the initial
seed of $\Fl_{w_0\times\nu}^{w_0\times(\lambda+\nu)}$
(already in the simplest example of 2-dimensional slices for
$\fg=\mathfrak{sl}_2$ where both variables are frozen,
cf.~\cite[Example~2.12]{kwy}).

We consider the following action of $\BZ^I$ on $\dZ^\alpha$: the generator
$(0,\ldots,0,1,0,\ldots,0)$ (1 at the $i$-th place) acts in the Pl\"ucker
coordinates $(Q_j,R_j)_{j\in I}$ of~\refss{morp} by an automorphism
$\eta_i(Q_j,R_j)=(Q_j,z^{\delta_{ij}}R_j-\delta_{ij}r_{i,a_i-1}Q_j)$ (it is easy
to check that this is indeed a biregular automorphism of $\dZ^\alpha$).
The frozen variables of the cluster structure on $\dZ^\alpha$ are
$(F_\alpha,q_{j,0})_{j\in I}$ where $F_\alpha$ is the equation of the boundary
of zastava, see e.g.~\cite[Section~5]{bdf}. Clearly, $\eta_i$ takes
$(F_\alpha,q_{j,0})_{j\in I}$ to $(F_\alpha q_{i,0}^{d_i},q_{j,0})_{j\in I}$, i.e.
does not preserve the frozen variables. However, it seems likely that
$\eta_i$ is an {\em almost cluster} transformation: a composition of a few
mutations, and the above change of frozen variables. Furthermore, if we set
$\eta_\nu:=\prod_{i\in I}\eta_i^{n_i}$ for $\nu=\sum_{i\in I}n_i\alpha_i$, then
the cluster structure transferred to $\dZ^\lambda$ from the isomorphism
$s^{\lambda^*+\nu^*}_{\nu^*}\circ\pr\colon \Fl_{w_0\times\nu^*}^{w_0\times(\lambda^*+\nu^*)}
\iso\dZ^\lambda$ differs from the reference one (transferred
from the isomorphism $s^{\lambda^*}_0\circ\pr\colon \Fl_{w_0}^{w_0\times\lambda^*}
\iso\dZ^\lambda$) by the automorphism $\eta_{\nu^*}$.

%Still, we hope that $\sigma_0^\nu$ takes
%the initial cluster of $\Fl_{w_0\times\nu}^{w_0\times(\lambda+\nu)}$ to some
%other cluster of $\Fl_{w_0}^{w_0\times\lambda}$, and thus all the clusters
%transferred to $\dZ^\alpha$ from the isomorphisms
%$s^{\lambda+\nu}_\nu\circ\pr\colon \Fl_{w_0\times\nu}^{w_0\times(\lambda+\nu)}
%\iso\dZ^\alpha$ are parts of a single cluster structure on $\dZ^\alpha$.

\ssec{sl2}{$\fg=\mathfrak{sl}_2$}
For $\fg=\mathfrak{sl}_2$, a positive coroot $\alpha$ is but a positive
integer $a$, and a cluster structure on $\dZ^a$ was defined
in~\cite[Section~5]{agsv} (where $\dZ^a$ is denoted $\CR_a$).
According to~\reft{sl2hank}, this cluster structure is a particular case of
the one of~\refss{desta}. In particular,
the exchange matrix $B(\varepsilon),\ (\varepsilon)=(2,0,\ldots,0)$
of~\cite[(5.16)]{agsv} coincides with the exchange matrix
of~\refss{exch}. Note that the cluster variables of~\cite{agsv} are certain
minors of a Hankel matrix composed of the coefficients of the formal series
$R(z)/Q(z)\in\BC[[z^{-1}]]$ (where $R,Q$ are the Pl\"ucker coordinates
of~\refss{morp}).
It would be nice to have such an explicit formula for
the cluster variables for general $\fg$. Also, the automorphism $\eta$
of~\refss{desta} is nothing but the {\em shift} of~\cite[Lemma~5.4.(i)]{agsv}
(a transformation from the type $A_{a-1}$ $Q$-system, cf. the paragraph
before~\cite[Remark~6.2]{agsv}).

\ssec{super}{Gaiotto-Witten superpotential}
%Recall the Cartan involution $\iota\colon \oZ^\alpha\to\oZ^\alpha$
%of~\cite[Section~4]{bdf}. Note that it preserves the open subvariety
%of trigonometric monopoles $\dZ^\alpha\subset\oZ^\alpha$ and induces the
%same named Cartan involution $\iota\colon \dZ^\alpha\to\dZ^\alpha$.
Let $K_i(z),\ i\in I$, be a collection of monic polynomials,
$K_i(z)=z^{l_i}+\kappa_{i,l_i-1}z^{l_i-1}+\ldots+\kappa_{i,0}$.
The data of $\{K_i(z)\}_{i\in I}$ is equivalent to the data of\\
(a) an ordered collection $\Lambda$ of dominant coweights
$\lambda_1,\ldots,\lambda_N$;\\
(b) an ordered configuration $(z_1,\ldots,z_N)$ of points in $\BA^1$.\\
Namely, given the above data we set
$K_i(z):=\prod_{1\leq n\leq N}(z-z_n)^{\langle\lambda_n,\check\alpha_i\rangle}$.
We denote by $\oA^\Lambda$ the moduli space of the above configurations of
{\em distinct} points $z_n$.

Recall the Gaiotto-Witten superpotential
$\CW^{\Lambda,\alpha}_-$: a multivalued holomorphic function on $\fh^\vee\times\oZ^\alpha\times\oA^\Lambda$
(see e.g.~\cite[1.8]{bdf}). We will denote by $\ol\CW{}^{\Lambda,\alpha}_-$
the restriction of $\CW^{\Lambda,\alpha}_-$
to $0\times\oZ^\alpha\times\oA^\Lambda$.
In the coordinates $w_{i,r},y_{i,r}$ of~\refss{morp} we have
\eq{gw}
\ol\CW{}^{\Lambda,\alpha}_-(\unl{w},\unl{y},\unl{K})=
\sum_{i,r}\frac{y_{i,r}K_i(w_{i,r})}{Q'_i(w_{i,r})}-\log F_\alpha+
\sum_{1\leq m<n\leq N}\lambda_m\cdot\lambda_n\log(z_m-z_n)
\end{equation}
where $F_\alpha$ is the equation of the zastava boundary
$\partial Z^\alpha=Z^\alpha\setminus\oZ^\alpha$ (see e.g.~\cite[Section~5]{bdf}).
Let $\BC[[z^{-1}]]\ni\frac{R_i(z)}{Q_i(z)}=:\sum_{p=0}^\infty h_{i,p}z^{-p-1}$.
Then
\eq{w}
\ol\CW{}^{\Lambda,\alpha}_-(\unl{w},\unl{y},\unl{K})=
\sum_{i,p}h_{i,p}\kappa_{i,p}-\log F_\alpha+
\sum_{1\leq m<n\leq N}\lambda_m\cdot\lambda_n\log(z_m-z_n)
\end{equation}

In case $\fg=\mathfrak{sl}_2$, the boundary equation $F_a$ is a frozen
cluster variable of the cluster structure on $\dZ^a$ of~\cite[Section~5]{agsv},
and all the coefficients $h_p$ are cluster variables according
to~\cite[Lemma~5.3,~Proposition~5.4]{agsv}. Hence
$\ol\CW{}^{\Lambda,\alpha}_-|_{K=K_0}$ is a constant
$\ell(K_0):=\sum_{1\leq m<n\leq N}\lambda_m\cdot\lambda_n\log(z_{0,m}-z_{0,n})$
plus a totally
positive function on $\dZ^a$ for a monic polynomial $K_0$ with nonnegative
coefficients $\kappa_p$. We expect a similar positivity for general $\fg$
(in particular, $F_\alpha$ is a frozen cluster variable). It would be 
interesting to study its tropicalization and the corresponding set of
positive integral tropical points, cf.~\cite{gs}.

\sec{append}{Appendix. $G=SL_2$: identification with the cluster structure
of~\cite{agsv}}

\centerline{\footnotesize GALYNA DOBROVOLSKA}

\medskip

Recall that $\slhat_2$ has two fundamental representations, which we denote by $V_{\varpi_1}$ and $V_{\varpi_0}$ in accordance with the notation of ~\refss{gaff}. Recall the generalized minors
$\varDelta_{w_0s_{i_1}\ldots s_{i_r}\varpi_{i_r},w_0\varpi_{i_r}}$ from ~\refss{seed}.
%which arose there as the initial cluster variables from \cite{l}.
Since for $\fg={\mathfrak{sl}}_2$ we have $w_0=s_1$, these generalized minors are of the form $\varDelta_{s_1(s_0 s_1)^r \varpi_1,s_1\varpi_1}$ and
$\varDelta_{s_1(s_0 s_1)^r\varpi_0,s_1\varpi_0}$.

%According to ~\refe{vardel}, ~\refe{vardels}, and ~\refp{comp}, we have that
%for an element $g(z)\in\widehat{SL}_2$, the functions $Q=Q_1$ and $R=R_1$
%introduced in ~\refss{morp} (evaluated at the point of the Zastava space
%which is the image of $g(z)$ under the embedding of ~\refp{coinc}) coincide
%with its lower right and lower left matrix coefficients, respectively.

Given a pair of polynomials $Q(z)=z^a +q_{a-1}z^{a-1} + ... + q_0,\
R(z)=r_{a-1}z^{a-1} + r_{a-2}z^{a-2} + ... + r_0$ representing a point of
$\dZ^a$, we can find a unique pair of polynomials
$F(z)=z^a +f_{a-1}z^{a-1} + ... + f_0,\
D(z)=d_{a-1}z^{a-1} + ... + d_1z + d_0$ such that $QF-RD=z^{2a}$.
Then both the matrix $g(Q,R):=\left(\begin{array}{cc}
z^{-a}F(z) & z^{-a}D(z)\\
z^{-a}R(z) & z^{-a}Q(z)
\end{array}\right)\in\widehat{SL}_2$
and its inverse matrix $g(Q,R)^{-1}=\left(\begin{array}{cc}
z^{-a}Q(z) & - z^{-a}D(z)\\
- z^{-a}R(z) & z^{-a}F(z)
\end{array}\right)\in\widehat{SL}_2$ actually lie in
$\bN_-\subset\widehat{SL}_2$
(notations of~\refss{gaff}). Moreover, according to~\refp{comp}, we have
$s^\lambda_0\circ\on{pr}(\overline{g(Q,R)})=(Q,R)$. Here $\lambda=a\alpha$
is a multiple of the simple coroot of ${\mathfrak{sl}}_2$, and
$\overline{g(Q,R)}\in\Fl^{s_1\times a\alpha}_{s_1}\subset\bFl_{s_1}=\bN_-\cdot s_1$
is the point $g(Q,R)\cdot s_1$.

On the other hand, following~\cite{agsv},
we consider the Taylor expansion at $\infty\in\BP^1$ of $\frac{R(z)}{Q(z)}=\frac{c_0}{z}+\frac{c_1}{z^2}+...+\frac{c_j}{z^j}+...$. We form the corresponding Hankel matrix using the elements $c_0,...,c_{2a-2}$, namely the $a \times a$ matrix $[c_{j+k}]_{j,k=0}^{a-1}$. We consider two kinds of minors of this matrix, the principal minors $\sC_1,...,\sC_a$ of sizes $1,...,a$, and the minors $\sD_1,...,\sD_{a-1}$ of sizes $1,...,a-1$ which are obtained from the principal minors of the same size by a shift of all entries by one unit to the right (or, equivalently, by a shift of all entries by one unit down). We will also denote these minors by
$\sC_r(Q,R),\ \sD_r(Q,R)$ when we want to stress the dependence on $Q,R$.
These Hankel minors (also called Hankel determinants) arise as cluster
variables in the cluster corresponding to $(\varepsilon)=(2,0,\ldots,0)$
in~\cite[5.2]{agsv}. See also the survey~\cite{hf} for general properties of Hankel matrices.

In this appendix we prove the following theorem:

\th{sl2hank} (a) The generalized minor
$\varDelta_{s_1(s_0 s_1)^r\varpi_1,s_1\varpi_1}(g(Q,R))$ is equal
(up to a change of sign) to the Hankel minor $\sC_r(Q,R)$.

(b) The generalized minor $\varDelta_{s_1(s_0 s_1)^r\varpi_0,s_1\varpi_0}(g(Q,R))$ is
equal (up to a change of sign) to the Hankel minor $\sD_r(Q,R)$.
\eth

Before starting the proof of~\reft{sl2hank}, we will recall a theorem of Kronecker which we will use.

First, for two polynomials $Q(z)=z^a +q_{a-1}z^{a-1} + ... + q_0,\
R(z)=r_{a-1}z^{a-1} + r_{a-2}z^{a-2} + ... + r_0$ we will write the $(2a-1)\times(2a-1)$ Sylvester matrix (the determinant of which computes the resultant of $Q$ and $R$) in the following form:
$$
\begin{bmatrix}
1 & q_{a-1} & \dots & \dots & q_0 & 0 & 0 & \dots & 0 & 0 \\
0 & 1 & q_{a-1} & \dots & q_1 & q_0 & 0 & \dots & 0 & 0 \\
\vdots & & \ddots & & & & & & \ddots & \\
0 & 0 & \dots & 1 & q_{a-1} & \dots & & \dots & q_1 & q_0 \\
0 & 0 & \dots & 0 & r_{a-1} & r_{a-2} & \dots & \dots & r_1 & r_0 \\
0 & 0 & \dots & r_{a-1} & r_{a-2} & \dots & & \dots & r_0 & 0 \\
\vdots & & & \dots & & & & & & \vdots \\
0 & r_{a-1} & \dots & \dots & r_0 &  0 & \dots & & \dots & 0\\
r_{a-1} & r_{a-2} & \dots & r_0 & 0 & 0 & \dots & & \dots & 0 \\
\end{bmatrix}
$$

Next we define for each $i$ an odd sub-resultant $\sR_i$ (which coincides with the notion of sub-resultant in~\cite[(2.6)]{uc}) to be the minor of the Sylvester matrix which is obtained by removing the same number $i$ of columns and rows at the top, the bottom, the right, and the left. We also define an even sub-resultant $\sS_i$ to be the minor of the Sylvester matrix obtained by removing the middle row, the same number $i$ of rows at the top and the bottom, and removing $i$ columns at the left, and $i+1$ columns at the right.

Now we can state the following formula of Kronecker~\cite{k}
(cf.~\cite[Corollary~3.2]{uc} for a modern reference) expressing
sub-resultants in terms of Hankel determinants of the Taylor expansion of the ratio of two polynomials:

\prop{Kronecker} (L.~Kronecker) $\sR_i=\sC_{a-i}$.
\eprop

We will also recall some facts from the theory of infinite-dimensional Lie algebras which we will need in the course of the proof (we will follow the exposition in \cite{kr} and use the notation of {\em loc. cit.}).

We start with an infinite-dimensional vector space vector space $V= \oplus_j \mathbb C v_j$. For each $m \in \ZZ$ we have the infinite-dimensional vector space $F^{(m)}$ with a vacuum vector $\psi_m=v_m \wedge v_{m-1} \wedge ...$ and a basis given by $v_{i_0} \wedge v_{i_{-1}} \wedge ...$ (such that $i_0 > i_{-1} > ...$ and $i_k = k+m$ for $k\ll0$). We define a group $GL_\infty$ as the group of invertible infinite matrices with entries $a_{i,j}$ ($i,j \in \mathbb Z$) such that all but finitely many of $a_{i,j}-\delta_{i,j}$ are zero. The group $GL_\infty$ acts
in $F^{(m)}$ as follows: $A(v_{i_0} \wedge v_{i_{-1}} \wedge ...)= A v_{i_0} \wedge A v_{i_{-1}} \wedge ... = \sum_{j_0 > j_{-1} > ...}
\det\!\! A_{j_0, j_{-1},...}^{i_0,i_{-1},...} v_{j_0} \wedge v_{j_{-1}} \wedge ...$, where $A_{j_0, j_{-1},...}^{i_0,i_{-1},...}$ denotes the matrix located on the intersection of the rows $j_0,j_{-1},...$ and columns $i_0,i_{-1},...$ of the matrix $A$.

Consider the standard $n$-dimensional representation $U$ of ${\mathfrak{sl}}_n$
with basis $u_1, u_2, ..., u_n$ (in this appendix we will only use $n=2$). Note that according to~\cite[(9.8)]{kr}, one can define an action of
$\slhat_n$ on $V$ in the following way, where an element
$A=\sum_i A_i z^i \in\slhat_n$ acts as the infinite matrix below (note that this action is obtained from the representation $U[z^{\pm1}]$ of $\slhat_n$ by identifying a basis element $v_j$ of $V$ with $z^k \cdot u_r$ where $j=kn+r$,
and $r \in \{1,2,...,n\}$):
$$
\begin{bmatrix}
 \dots & \dots & \dots & \dots & \dots & \dots \\
\dots & A_{-1} & A_0 & A_{1} & \dots & \dots \\
\dots & \dots & A_{-1}& A_0 & A_1 & \dots \\
\dots & \dots & \dots & A_{-1}& A_0 & \dots \\
 \dots & \dots & \dots & \dots & \dots & \dots \\
\end{bmatrix}
$$

Note that the images of the matrices in $\slhat_n$ obtained in this way have finitely many non-zero diagonals. Therefore by \cite[Section~4.4]{kr}
the action of $\slhat_n$ in $F^{(m)}$ is given by the same formula as for $GL_\infty$, namely $A(v_{i_0} \wedge v_{i_{-1}} \wedge ...) =  \sum_{j_0 > j_{-1} > ...}
\det\!\! A_{j_0, j_{-1},...}^{i_0,i_{-1},...} v_{j_0} \wedge v_{j_{-1}} \wedge ...$ for $A \in\slhat_n$.
This way for $m=0,1,2,...,n-1$ we obtain all the fundamental representations $V_{\varpi_m}$, where $V_{\varpi_m}$ is the irreducible sub-representation of $\slhat_n$ in $F^{(m)}$ which is generated by the vacuum vector $\psi_m=v_m \wedge v_{m-1} \wedge ...$.

Finally, note that the action of $\mathfrak n_{-}$ in $F^{(m)}$ is not integrable (in general, $g v_{i_0} \wedge v_{i_{-1}} \wedge ...$ is an infinite sum for an element $g \in \mathbf{N}_{-}$). However, for any basis element $v_{j_0} \wedge v_{j_{-1}} \wedge ...$ of $F^{(m)}$ its coefficient in the (infinite) expansion of $g v_{i_0} \wedge v_{i_{-1}} \wedge ...$ in the elements of the basis of $F^{(m)}$ is well-defined and can be computed by the same formula as for $g \in GL_\infty$, namely it is equal to $\det\!\! A_{j_0, j_{-1},...}^{i_0,i_{-1},...}$. Note that here we caculate $\det\!\! A_{j_0, j_{-1},...}^{i_0,i_{-1},...}$ in the following way. By definition there exists $N$ such that for $k < N$ we have $j_k=i_k=k+m$; then $\det\!\! A_{j_0, j_{-1},...}^{i_0,i_{-1},...} =  \det\!\! A_{j_0, j_{-1},..., j_N}^{i_0,i_{-1},...,j_N}$, which is a finite determinant. The justification for this formula is that the infinite matrix with rows $j_0, j_{-1},...$ and columns $i_0, i_{-1},...$ can be divided into four blocks, where the two diagonal blocks are the finite block with rows $j_0, j_{-1},...,j_N$ and columns $i_0, i_{-1},...,i_N$ and the inifinite lower-triangular block with $1$'s on the diagonal with rows $j_{N-1}, j_{N-2}, ...$ and columns $i_{N-1}, i_{N-2}, ...$.

%According to~\cite[Section~6.2]{kr}, a larger group $\overline{GL}_\infty$ (consisting of of infinite matrices with entries $a_{i,j}$ %such that all but finitely many $a_{i,j}-\delta_{i,j}$ with $i\geq j$ are zero) also acts on $F^{(m)}$ by the same formula %$R(A)%(v_{i_0} \wedge v_{i_{-1}} \wedge ...) = \sum_{j_0 > j_{-1} > ...} \mathrm d \mathrm e \mathrm t A_{j_0, %j_{-1},...}^{i_0,i_{-1},...} v_{j_0} \wedge v_{j_{-1}} \wedge ...$.

Now we are ready to prove our theorem.

\begin{proof}[Proof of~\reft{sl2hank}]

Using exterior powers, the computation of generalized minors reduces to the computation of finite minors because the infinite matrices we use are upper triangular up to a finite portion. As a result we obtain finite minors which stand in certain rows and columns of the infinite matrix, depending on the element of the Weyl group which appears in the definition of a particular generalized minor.

For example, for the element of the Weyl group given by $s_0 s_1 s_0 s_1$ and the fundamental representation $V_{\varpi_1}$, we obtain the following minors:

\tiny

$$
\kbordermatrix{
 \ & -5 & -4 & {\bf -3} & {\bf -2} & {\bf -1} & {\bf 0} & 1 & {\bf 2} & 3 & 4 & 5 & 6 & 7 & 8\\
 -5 & 1  & \   & \   & \    & \   & \  & \  & \  & \ & \  & \ & \ & \  & \ \\
 -4 & 0  & 1   & \           & \            & \             & \          & \  & \           & \ & \  & \ & \  & \  & \  \\
 -3 & f_{a-1}  & d_{a-1}   & 1          & \            & \             & \          & \  & \           & \ & \  & \ & \  & \  & \  \\
 {\bf -2} & r_{a-1}  & q_{a-1}   & {\bf 0}          & {\bf 1}          & \             & \          & \  & \           & \ & \  & \ & \  & \  & \  \\
 -1 &f_{a-2}  & d_{a-2} & f_{a-1}  & d_{a-1}      & 1             & \          & \  & \           & \ & \  & \ & \  & \  & \  \\
  {\bf 0} & r_{a-2}  & q_{a-2}&\bf r_{a-1}  &\bf q_{a-1}   & \bf  0             & \bf 1          & \  & \           & \ & \  & \ & \  & \  & \  \\
  1 &f_{a-3} & d_{a-3}&f_{a-2}  & d_{a-2} &f_{a-1}  & d_{a-1} & 1  & \           & \ & \  & \ & \  & \  & \  \\
  {\bf 2} &r_{a-3}&q_{a-3}&\bf r_{a-2}  &\bf  q_{a-2}&\bf r_{a-1}  &\bf  q_{a-1}   & 0  & \bf 1           & \ & \  & \ & \  & \  & \  \\
 3 &f_{a-4}&d_{a-4}&f_{a-3} & d_{a-3}& f_{a-2}  & d_{a-2} &f_{a-1}  & d_{a-1} & 1   & \  & \ & \  & \  & \  \\
 {\bf 4} &r_{a-4}&q_{a-4}&\bf r_{a-3}&\bf q_{a-3}&\bf r_{a-2}  &\bf  q_{a-2}& r_{a-1}  &\bf q_{a-1}   & 0  & 1 & \ & \  & \  & \  \\
 5 &f_{a-5}&d_{a-5}&f_{a-4}&d_{a-4}&f_{a-3} & d_{a-3}&f_{a-2}  & d_{a-2} &f_{a-1}  & d_{a-1}& 1   & \  & \  & \  \\
 {\bf 6} &r_{a-5}&q_{a-5}&\bf r_{a-4}&\bf q_{a-4}&\bf r_{a-3}&\bf q_{a-3}&r_{a-2}&\bf q_{a-2}&r_{a-1}&q_{a-1}& 0&1&\ &\  \\
 7 &f_{a-6}&d_{a-6}&f_{a-5}&d_{a-5}&f_{a-4}&d_{a-4}&f_{a-3}&d_{a-3}&f_{a-2}&d_{a-2}&f_{a-1}&d_{a-1}&1& \  \\
 8 &r_{a-6}&q_{a-6}&r_{a-5}&q_{a-5}&r_{a-4}&q_{a-4}&r_{a-3}&q_{a-3}&r_{a-2}&q_{a-2}& r_{a-1}& q_{a-1}& 0&1  }
$$

\normalsize

After we collect the entries at the intersections of the marked rows and columns, we obtain the following submatrix (which after transposing it and exchanging the order and the signs of some the rows will be exactly the $(a-3)$-th sub-resultant $\sR_{a-3}$ for the polynomials $Q$ and $R$):

\tiny

$$
\begin{bmatrix}
0 & 1 & \ & \ & \ \\
r_{a-1} & q_{a-1} & 0 & 1 & \ \\
r_{a-2} & q_{a-2} & r_{a-1} & q_{a-1} & 1 \\
r_{a-3} & q_{a-3} & r_{a-2} & q_{a-2} & q_{a-1} \\
r_{a-4} & q_{a-4} & r_{a-3} & q_{a-3} & q_{a-2} \\
\end{bmatrix}
$$
\normalsize

We see that the finite minors we obtain up to permutation of rows and transposition are exactly the sub-resultants in the case of odd number of rows. Indeed, for the element of the Weyl group given by $(s_0 s_1)^r$ and the fundamental representation $V_{\varpi_1}$ we obtain that the finite minor consists of $2r+1$ rows numbered $2r+2,2r,\ldots,-2r+2$ and $2r+1$ columns numbered $2,0,-1,\ldots,-2r+2$. If we permute the columns so that the odd rows stand on the left in the same order as they were before and even rows stand on the right in the same order as they were before, we obtain exactly an odd sub-resultant as defined above. Now we use the above theorem of Kronecker to conclude that the generalized minors for the fundamental representation $V_{\varpi_1}$ of $\slhat_2$ are equal to the corresponding principal Hankel minors.

Now we turn to the generalized minors of the basic representation $V_{\varpi_0}$.
For example, for the element $s_0 s_1 s_0 s_1$ of the Weyl group and the fundamental representation $V_{\varpi_0}$ we have the following submatrix:

\tiny

$$
\kbordermatrix{
 \ & -5 & -4 & {\bf -3} & {\bf -2} & {\bf -1} & {\bf 0} & 1 & 2 & 3 & 4 & 5 & 6 & 7 & 8\\
 -5 & 1  & \   & \   & \    & \   & \  & \  & \  & \ & \  & \ & \ & \  & \ \\
 -4 & 0  & 1   & \           & \            & \             & \          & \  & \           & \ & \  & \ & \  & \  & \  \\
 -3 & f_{a-1}  & d_{a-1}   & 1          & \            & \             & \          & \  & \           & \ & \  & \ & \  & \  & \  \\
 {\bf -2} & r_{a-1}  & q_{a-1}   & {\bf 0}          & {\bf 1}          & \             & \          & \  & \           & \ & \  & \ & \  & \  & \  \\
 -1 &f_{a-2}  & d_{a-2} & f_{a-1}  & d_{a-1}      & 1             & \          & \  & \           & \ & \  & \ & \  & \  & \  \\
  {\bf 0} & r_{a-2}  & q_{a-2}&\bf r_{a-1}  &\bf q_{a-1}   & \bf  0             & \bf 1          & \  & \           & \ & \  & \ & \  & \  & \  \\
  1 &f_{a-3} & d_{a-3}&f_{a-2}  & d_{a-2} &f_{a-1}  & d_{a-1} & 1  & \           & \ & \  & \ & \  & \  & \  \\
  {\bf 2} &r_{a-3}&q_{a-3}&\bf r_{a-2}  &\bf  q_{a-2}&\bf r_{a-1}  &\bf  q_{a-1}   & 0  & 1           & \ & \  & \ & \  & \  & \  \\
 3 &f_{a-4}&d_{a-4}&f_{a-3} & d_{a-3}& f_{a-2}  & d_{a-2} &f_{a-1}  & d_{a-1} & 1   & \  & \ & \  & \  & \  \\
 {\bf 4} &r_{a-4}&q_{a-4}&\bf r_{a-3}&\bf q_{a-3}&\bf r_{a-2}  &\bf  q_{a-2}& r_{a-1}  & q_{a-1}   & 0  & 1 & \ & \  & \  & \  \\
 5 &f_{a-5}&d_{a-5}&f_{a-4}&d_{a-4}&f_{a-3} & d_{a-3}&f_{a-2}  & d_{a-2} &f_{a-1}  & d_{a-1}& 1   & \  & \  & \  \\
 6 &r_{a-5}&q_{a-5}&r_{a-4}&q_{a-4}&r_{a-3}&q_{a-3}&r_{a-2}&q_{a-2}&r_{a-1}&q_{a-1}& 0&1&\ &\  \\
 7 &f_{a-6}&d_{a-6}&f_{a-5}&d_{a-5}&f_{a-4}&d_{a-4}&f_{a-3}&d_{a-3}&f_{a-2}&d_{a-2}&f_{a-1}&d_{a-1}&1& \  \\
 8 &r_{a-6}&q_{a-6}&r_{a-5}&q_{a-5}&r_{a-4}&q_{a-4}&r_{a-3}&q_{a-3}&r_{a-2}&q_{a-2}& r_{a-1}& q_{a-1}& 0&1  }
$$

\normalsize

After we collect the entries at the intersections of the marked rows and columns we obtain the following submatrix (which after transposing and exchanging the order of the rows is the $(a-3)$-th sub-resultant $\sS_{a-3}$ for the polynomials $Q$ and $R$):

\tiny

$$
\begin{bmatrix}
0 & 1 & \ & \  \\
r_{a-1} & q_{a-1} & 0 & 1  \\
r_{a-2} & q_{a-2} & r_{a-1} & q_{a-1}  \\
r_{a-3} & q_{a-3} & r_{a-2} & q_{a-2}  \\
r_{a-4} & q_{a-4} & r_{a-3} & q_{a-3} \\
\end{bmatrix}
$$
\normalsize

More generally, for the element of the Weyl group given by $(s_0 s_1)^r$ and the fundamental representation $V_{\varpi_0}$ we obtain that the finite minor consists of $2r+1$ rows numbered
$2r,2r-2,\ldots,-2r+2$ and $2r$ rows numbered $0,-1,\ldots,-2r+1$. If we permute the columns so that the odd rows stand on the left in the same order as they were before and even rows stand on the right in the same order as they were before, we obtain exactly an even sub-resultant as defined above.

Finally, we can reduce the case of even number of rows to the case of odd number of rows in the following way. Note that the even sub-resultant $\sS_i$ of the polynomials $R(z)$ and $Q(z)$ is equal to a usual sub-resultant
(cf.~\cite[(2.6)]{uc}) of the polynomials $R(z)$ and $z Q(z)$. We claim that this (usual) sub-resultant of the polynomials $R(z)$ and $z Q(z)$ is equal to the determinant of the original Hankel minor $D_{a-i-1}$. Indeed, by the theorem of Kronecker mentioned above, the usual sub-resultants of $R(z)$ and $z Q(z)$ equal the corresponding (principal) Hankel minors for $R(z)$ and $z Q(z)$. To conclude, we notice that the equality $R(z)/Q(z)=z(R(z)/(z Q(z)))$ implies that the principal Hankel minors for $R(z), z Q(z)$ are equal to the Hankel minors for $R(z),Q(z)$ which are obtained by shifting the corresponding principal Hankel minors by one unit to the right.
\end{proof}

\bigskip

\footnotesize{{\bf M.F.}: National Research University Higher
School of Economics, Russian Federation,\\
Department of Mathematics, 6 Usacheva st, Moscow 119048;\\
Skolkovo Institute of Science and Technology;\\ %\hphantom{x}\ab\,
{\tt fnklberg@gmail.com}}

\footnotesize{
{\bf A.K.}: Steklov Mathematical Institute, Algebraic Geometry Section,
8 Gubkina, Moscow 119991, Russia;\\
Laboratory of Algebraic Geometry,\\
National Research University Higher
School of Economics, Russian Federation;\\
{\tt akuznet@mi.ras.ru}}

\footnotesize{{\bf L.R.}: National Research University Higher
School of Economics, Russian Federation,\\
Department of Mathematics, 6 Usacheva st, Moscow 119048;\\
Institute for the Information Transmission Problems of RAS;\\ %\hphantom{x}\ab\,
{\tt leo.rybnikov@gmail.com}}

\footnotesize{
{\bf G.D.}: Department of Mathematics, Columbia University, \\
2990 Broadway, New York, NY 10027, USA; \\
{\tt galdobr@gmail.com}}

\end{document}